\documentclass[12pt]{amsart}
\usepackage{amsmath,amsthm,amsfonts,amssymb}
\usepackage{amsmath, amsthm, amssymb, amscd, amsfonts}
\usepackage[all]{xy}
\usepackage{amssymb, amsthm, amsmath}
\usepackage[english]{babel}
\usepackage{amsfonts}
\usepackage[dvipsnames]{xcolor}
\usepackage[colorlinks=true, pdfstartview=FitV, linkcolor=blue, 
            citecolor=blue, urlcolor=blue]{hyperref}
\newcommand{\bab}{\color{DarkOrchid}{}}
\newcommand{\eab}{\normalcolor{}}
\newcommand{\ssha}{\;semisimple Hopf agebra\; }
\newcommand{\thn}{Then\;}
\newcommand{\stat}{\mdn {\bf \huge Statement}\;}
\usepackage{color}
\DeclareMathOperator{\ad}{ad} 
\DeclareMathOperator{\id}{id}
\DeclareMathOperator{\Hom}{Hom}

\newcommand{\ra}{\rightarrow}
 \newcommand{\ovr}{\overline}

\newcommand{\Z}{\mathbb Z}

\newcommand{\ot}{\otimes}
\newcommand{\mtc}{\mathcal}
\newcommand{\lam}{\lambda}
\newcommand{\lb}{\label}
\newcommand{\Lam}{\Lambda}

\newcommand{\al}{\alpha}
\newcommand{\onh}{On the other hand\;}
\newcommand{\eps}{\epsilon}

\newcommand{\D}{\Delta}
\newcommand{\ro}{\rho}
\newcommand{\teta}{\theta}

\newcommand{\rh}{\rightharpoonup}
\newcommand{\lh}{\leftharpoonup}
\numberwithin{equation}{section}
\newtheorem{lem}[equation]{Lemma}

\newtheorem{thm}[equation]{Theorem}
\newtheorem{prop}[equation]{Proposition}
\newtheorem{defn}[equation]{Definition}
\newtheorem{cor}[equation]{Corollary}
\newtheorem{rem}[equation]{Remark}

\newcommand{\dw}{\downarrow}
\newcommand{\uw}{\uparrow}

\newcommand{\ch}{\chi}
\newcommand{\mtr}{\mathrm}
\newcommand{\co}{\mathcal{O}}

\numberwithin{equation}{section}
\newcommand{\ncm}{\newcommand}
\ncm{\np}{\newpage}
\ncm{\ebl}{\end{thebibliography}}
\ncm{\bbl}{\begin{thebibliography}}
\ncm{\chd}{_{ _{\ch}}}
\ncm{\ald}{_{ _{\al}}}
\ncm{\cP}{\mathcal{P}}
\ncm{\ei}{e_i}
\ncm{\eij}{e_{i,\;j}}
\newcommand{\bb}{\blue}
\ncm{\bt}{\begin{thm}}
\ncm{\bdef}{\begin{defn}}
\ncm{\edf}{\end{defn}}
\ncm{\et}{\end{thm}}
\ncm{\bc}{\begin{cor}}
\ncm{\bl}{\begin{lem}}
\ncm{\el}{\end{lem}}
\ncm{\bpf}{\begin{proof}}
\ncm{\epf}{\end{proof}}
\ncm{\ec}{\end{cor}}
\ncm{\er}{\end{rem}}
\ncm{\br}{\begin{rem}}
\ncm{\bn}{\begin}
\ncm{\ca}{\mtc{A}}
\ncm{\cc}{\mtc{C}}
\ncm{\cz}{\mtc{Z}}
\newcommand{\blu}{\blue}
\ncm{\bp}{\begin{prop}}
\ncm{\ep}{\end{prop}}
\ncm{\bd}{
\begin{document}}
\ncm{\ed}{\end{document}}
\ncm{\beq}{\begin{equation}}
\ncm{\beqn}{\begin{equation*}}
\ncm{\eeq}{\end{equation}}
\ncm{\eeqn}{\end{equation*}}
\ncm{\bea}{\begin{eqnarray}}
\ncm{\eea}{\end{eqnarray}}
\ncm{\beanon}{\begin{eqnarray*}}
\ncm{\eeanon}{\end{eqnarray*}}\ncm{\ek}{\eps|_K}\ncm{\diez}{\#}
\newcommand{\ck}{\mathcal{K}}
\ncm{\bwt}{\bowtie}
\ncm{\cC}{\mtc{C}}
\ncm{\cX}{\mtc{X}}
\newcommand{\K}{\mathrm{K}}
\ncm{\wt}{\widetilde}
\ncm{\sg}{\sigma}\ncm{\Rep}{\mathrm{Rep}}
\ncm{\Irr}{\mathrm{Irr}}\ncm{\X}{\mathcal{X}}
\ncm{\cA}{\mathcal{A}}
\ncm{\HKer}{\mtr{HKer}}
\ncm{\LKER}{\mtr{LKer}}
\ncm{\aad}{\mtr{ad}}
\ncm{\Dr}{\mtr{D}}
\ncm{\cD}{\mathcal{D}}
\ncm{\G}{\mathcal{G}}
\ncm{\Dc}{\mtc{D}}
\ncm{\E}{\mtc{E}}
\ncm{\fp}{\mtr{FPdim}}
\ncm{\Vc}{\mtr{Vec}}
\ncm{\cK}{\mtc{K}}
\ncm{\cM}{\mtc{M}}
\newcommand{\bq}{\beq}
\newcommand{\eq}{\eeq}
\ncm{\cE}{\mtc{E}}
\ncm{\cS}{\mtc{S}}
\ncm{\End}{\mtr{End}}
\ncm{\cop}{\mtr{cop}}
\ncm{\op}{\mtr{op}}
\ncm{\chr}{character }
\ncm{\bw}{\bwt}
\newcommand{\core}{\mtr{core}}
\newcommand{\ce}{\cE}
%\ncm{\Vec}{\mtr{Vec}}
\ncm{\cY}{\mtc{Y}}
\newcommand{\cx}{\mathcal{X}}
\ncm{\hker}{\mtr{HKer}}
\ncm{\bx}{\boxtimes}\ncm{\cd}{\cD}
\ncm{\blue}{\textcolor[rgb]{.00, .00, 1.00}}
\ncm{\red}{\textcolor[rgb]{1.00, .00, .00}}
\ncm{\green}{\textcolor[rgb]{.00, 1.00, .00}}
\ncm{\bne}{\begin{enumerate}}
\ncm{\ene}{\end{enumerate}}
\ncm{\lker}{\mtr{LKer}}
\newcommand{\irr}{\Irr}
\ncm{\md}{\medbreak}
\ncm{\rep}{\Rep} \ncm{\cs}{\mtc{S}}
\newcommand{\tet}{\theta}
\newcommand{\kk}{\Bbbk}
\newcommand{\kt}{$\Bbbk$\nobreakdash-\hspace{0pt}}
\newcommand{\trait}{\nobreakdash-\hspace{0pt}}
\newcommand{\Rt}{$\mathrm{R}$\nobreakdash-\hspace{0pt}}
\newcommand{\ti}{\mbox{-}\,}
\newcommand{\N}{\mathbb{N}}
\newcommand{\mdn}{\medbreak\noindent}
\ncm{\cO}{\mtc{O}}
\newcommand{\xra}{\xrightarrow}

\title[Semisimple Hopf algebras]
{On a symmetry of  M\"uger's centralizer for the %representations of 
Drinfeld double of a semisimple Hopf algebra}
\author{Sebastian Burciu}

\address{Inst.\ of Math.\ ``Simion Stoilow" of the Romanian Academy, Research Unit 5\\ P.O. Box 1-764, RO-014700, Bucharest, Romania}
%\thanks{This work was supported by a grant of the Romanian National Authority for Scientific Research, CNCS Ð UEFISCDI, project number PN-II-RU-TE-2012-3-0168.}
\email{sebastian.burciu@imar.ro}
\date{\today}
\bd
\maketitle
\bn{abstract}
%Normal coideal subalgebras of semisimple Drinfeld doubles are studied. 
In this paper we prove a formula that relates M\"uger's centralizer in the category of representations of a factorizable Hopf algebra to the notion of Hopf kernel of a representation of the dual Hopf algebra.
 Using this relation we obtain a complete description for M\"uger's centralizer of some fusion subcategories of the fusion category  of finite dimensional representations of a  Drinfeld double of a semisimple Hopf algebra.
\end{abstract}
%\frontmatter
%\tableofcontents
\section{Introduction and main results}
%\subsection{Importance of M\"uger's centralizer}
M\"uger has introduced in \cite{proclond} the notion of centralizer of a fusion subcategory of a braided fusion category. One of the most remarkable features of this notion is that the centralizer of a nondegenerate fusion subcategory of a modular category is a categorical complement of the nondegenerate subcategory. This principle is applied in many classification results of fusion categories, see for example \cite{DGNO, DGNOI, ENO2}. 
\md
Despite its importance, in general it is a difficult task to give a concrete description for the centralizer of all fusion subcategories of a given fusion category. Only few cases are know in the literature.  For instance, in the same aforementioned  paper, \cite{proclond},  M\"uger described  the centralizer of all fusion subcategories of the category of finite dimensional representations of a Drinfeld double of a finite abelian group. More generally, for the category of representations of a (twisted) Drinfeld double of an arbitrary finite group a similar formula was then given in \cite{nnw}. For the braided center of Tambara-Yamagami categories, this centralizer was described by computing completely the $S$-matrix of the modular category in \cite{gnn}.
\md In this paper we study some properties of M\"uger's centralizer for the category of representations of a semisimple factorizable Hopf algebra. A formula that relates M\"uger's centralizer  to the Hopf kernel of representations of the dual Hopf algebra is proven in Theorem \ref{wding}. Then we specialize  these results to the category of representations of Drinfeld doubles of semisimple Hopf algebras. %This relation is based on the inequality for evaluation of characters at dual characters from \cite{Bker}.
This allows us to obtain a certain symmetry of M\"uger's centralizer for these categories. This symmetry can also be viewed as a  generalization of the above mentioned results for the (twisted) Drinfeld double of a  finite group.%in \cite{nnw}.}
\md
Normal Hopf subalgebras of a semisimple Drinfeld double were studied in \cite{nrm}. In loc.cit. the author shows that if $K, L$ are two normal Hopf subalgebras of semisimple Hopf algebra $A$ such that $[K, L]=0$ and $[(A//K)^{*}, (A//L)^{*}]=0$ then $B(K,L):=(A//K)^{*\cop}\bwt L$ is a normal Hopf subalgebra of $D(A)$. In this paper we denote by $\cd(K,L):=\rep(D(A)//B(K,L))$ the fusion subcategory of $\rep(D(A))$ obtained by taking the normal quotient $B(K,L)$.
\mdn
Our first main result is the following symmetry for the centralizer of $\cd(K, L)$:
\bt\label{centraliz} Let $A$ be a semisimple Hopf algebra and $K,L$ be two normal Hopf subalgebras of $A$ such that $B(K,L)$ is a normal Hopf subalgebra of $D(A)$. Then
\beqn
\cd(K, L)'=\cd(L,K)
\eeqn
as fusion subcategories of $\rep(D(A))$.
\et
Let $A$ be a semisimple Hopf algebra and $K$ be a Hopf subalgebra of  the Drinfeld double $D(A)$. We denote by $\cd(K)$ the fusion subcategory of $\rep(D(A))$ whose objects  are those $D(A)$-representations that receive a trivial $K$-action.% from the Hopf subalgebra $K$ of $D(A)$. 

Our second main result gives a description for the centralizer of $\cd(K)$:
\bt\label{justonecomp}
Let $A$ be a semisimple Hopf subalgebra and $K$ be a  normal Hopf subalgebra of $A$. Then 
\beqn
\cD(K)'=<K>
\eeqn
where $<K>$ is the fusion subcategory of $\Rep (D(A))$ generated by the $D(A)$-module $K$. Here $K$ is regarded as $D(A)$-module via the action
\beq
(f \bwt a)x=a_{1}xS(a_{2})\lh S^{-1}f
\eeq
where $a \in A$, $f \in A^{*}$ and $x \in K$.
\et
%\blue{The above theorem follows form the following description for the centralizer of the fusion subcategory $\cd(K)$:
%See subsection \ref{mod} for more details on the structure of $K$ as a $D(A)$-module.}
 % \bb{As an application it shows that $\rep(D(A))$ contains a normal Lagrangian fusion subcategory if and only if $A$ is an abelian extension.}
\subsection*{Organization of the paper}In Section \ref{prelim} we recall the definition of M\"uger's centralizer of a subcategory of a modular fusion category as well as its basic properties that are used throughout this paper. Section \ref{lattice} presents some results concerning the lattice of fusion subcategories of the category of representation of a semisimple Hopf algebra. We recall here Brauer's theorem for kernels of representations, as presented in \cite{gmj}.  In Section \ref{kcc} we prove Theorem \ref{wding} that relates the kernel of a corepresentation and M\"uger's centralizer in the category of representations of a factorizable semisimple Hopf algebra.  Theorem \ref{justonecomp} and Theorem \ref{centraliz} for are proven in Section \ref{kccd}.% As an example, representations of finite groups.\md
\md
We work over an algebraically closed field $\kk$ of characteristic zero. We use Sweedler's notation for comultiplication with the sigma symbol dropped and the notation $S$ for the antipode. All the other Hopf algebra notations are those used in \cite{Montg}.
\subsection*{Acknowledgments}
This work was supported by a grant of the Romanian National Authority for Scientific Research, CNCS-UEFISCDI, project number PN-II-RU-TE-2012-3-0168.
\section{Preliminaries}\lb{prelim}

In this section we recall some preliminary results that are further needed throughout the paper.

\subsection*{General conventions on fusion categories} As usually, by a fusion category we mean a $\kk$-linear semisimple rigid tensor category $\cc$ with finitely many isomorphism
classes of simple objects, finite dimensional spaces of morphisms, and such that the unit object of $\cc$ is simple. We refer the reader to \cite{ENO} for a general theory of such
categories. 
\mdn For a fusion category $\mtc{C}$ we denote by $\Lam_{\cc}$ the set of isomorphism classes
of simple objects of $\cc$ and by $\co(\cc)$ the class of all objects of $\cc$.  Recall that the Grothendieck ring $K_0(\cc)$
of $\cc$ is the free $\Z$-module generated by 
$\Lam_{\cc}$ with the multiplication induced by the tensor product in $\cc$. The Grothendieck
ring $K_0(\cc)$ is a based unital ring  (see for example \cite{NG} for definition of based
rings). %The isomorphism classes $[X]$ of simple objects of $\cc$ form a $\Z$-basis for$K_0(\cc)$.
%A fusion category is called pointed if all its simple objects are invertible.
\mdn
Recall \cite{ENO} that there is a unique algebra homomorphism  $\fp: \mtr{K}_{0}(\cc)\ra \kk^{*}$ such that $\fp([X])>0$ for any simple object $X$ of $\cc$. By definition, the Frobenius-Perron dimension of $\cc$ is given by $\fp(\cc):=\sum_{X \in \Lam_{\cc}}\fp(X)^{2}$.
\mdn
By a fusion subcategory of a fusion category $
\cc$ we understand a full replete tensor subcategory
of $\cc$.  For two fusion subcategories $\cd$ and $\ce$ of $\cc$ we denote by $\cd \vee \ce $ the smallest fusion subcategory of $\cc$ containing both $\cd$ and $\ce$ as fusion subcategories. If $\ce\subseteq \cd$ are fusion subcategories of a fusion category $\cc$ such that $\fp(\cd)=\fp(\ce)$ then clearly $\cd=\ce$.
\subsection{Modular fusion categories} Recall that a braided tensor category $\cC$ is a tensor category equipped for all $X,Y \in \co(\cc)$ with natural isomorphisms $c_{X,\;Y} : X \otimes Y \ra  Y \otimes X$  satisfying the hexagon axiom, see for example \cite{BaKi, JS}. 
%\subsubsection*{Ribbon fusion categories}
\mdn
A twist on a braided fusion category $\cc$ is a natural automorphism $\tet:\id_{\cc}\ra \id_{\cc}$ satisfying $\tet_{1}=\id_{1}$
and
\beq
\tet_{X\ot Y}=(\tet_{X}\ot \tet_{Y})c_{YX}c_{XY}.
\eeq
A braided fusion category is called {\it premodular} or {\it ribbon} if it has a twist satisfying $\tet_{X^{*}} =\tet_{X}^{*}$ for all $X \in \co(\cc)$.
\mdn
Recall that the entries of the $S$-matrix,  $S=\{s_{X,Y}\}$ of a premodular category are defined as the quantum trace $s_{X,Y}:=tr_{q}(c_{YX}c_{XY})$, see \cite{tu}. A premodular category $\cc$ is called {\it  modular} if the above $S$-matrix is nondegenerate.

\subsection*{Centralizers in braided fusion categories}

Let $\cc$ be a modular fusion category and $\ck$ be a fusion subcategory of $\cc$. M\"uger introduced the notion of {\it centralizer of $\ck$ } as the fusion subcategory $\ck'$ of $\cc$ generated by all simple objects $X$ of $\cc$ satisfying  
\beq\lb{centr} c_{X,\; Y}c_{Y,\;X}=\mtr{id}_{X\otimes Y}\eeq
for all objects $Y \in \co(\ck)$ (see \cite{proclond}).
\mdn
In the case of a ribbon category $\cc$ the condition of Equation \eqref{centr} is equivalent to 
\beq\lb{centr2} s_{X,\;Y}=\fp(X)\fp(Y)\eeq for all objects $Y \in \co(\ck)$.
Note that  in  general \beq\label{is}|s_{X,\;Y}|\leq \fp(X)\fp(Y)\eeq  by \cite[Proposition 2.5]{proclond}. In the situation of Equation \eqref{centr} we say that the objects $X$ and $Y$ centralize each other.
\subsubsection*{\bf Properties of the centralizer}

If $\cC$ is a modular fusion category and $\cK$ a fusion subcategory of $\cC$ then by  \cite[Theorem 3.2]{proclond} one has $\cK''=\cK$ and \beq\fp(\cK)\fp(\cK')=\fp (\cC).\eeq Moreover by \cite[Theorem 4.2]{proclond} it follows that if $\cK$ is also a modular category then $\cC \cong \cK\boxtimes \cK'$.
Let $\cD$ and $\cE$ be fusion subcategories of a modular fusion category $\cC$. Then one has 
\beq\lb{centrofint} (\cD\vee\cE)'=\cD'\cap \cE' \text{\;and \;} (\cD\cap \cE)'=\cD'\vee \cE',\eeq see \cite{proclond}.
%\subsubsection{Properties of the $S$-matrix}
%$$s_{X,\;Y}=\teta^{-1}_X\teta^{-1}_Y\sum_{W\in \Lam_{\cc}}N^W_{XY}\teta_Wd(W)$$
%$$s_{X,\;Y} s_{X,\;Z} =d(X)\sum_{W\in \Lam_{\cc}}N^W_{YZ} s_{X,\;W},$$ for all $X,Y,Z\in \Lam_{\cc}$.

%\blue{Definition of the N's and references for any of the results. Clarify which ones are used.}
\ncm{\cZ}{\mathcal{Z}}
\ncm{\cL}{\mtc{L}}
%\subsection{Fusion subcategories of $\Rep(A)$}
\section{On the lattice of fusion subcategories of $\Rep(A)$}\label{lattice} Let $A$ be a finite dimensional semisimple Hopf
algebra over $\kk$. Then $A$ is also cosemisimple and $S^{2}=\id$, see \cite{Lard}. 
The character ring
$C(A)$ of $A$ is a semisimple $\kk$-subalgebra of $A^*$ and it has a vector space basis
given by the set $\mtr{Irr}(A)$ of irreducible characters of $A$, see  \cite{Z}. Moreover,
$C(A)=\mtr{Cocom}(A^*)$, the space of cocommutative elements of $A^*$. By duality, the
character ring of $A^*$ is a semisimple $\kk$-subalgebra of $A$ and $C(A^*)=\mtr{Cocom}(A)$. If
$M$ is an $A$-representation with character $\chi$ then $M^*$ is also an
$A$-representation with character $\chi^*=\chi \circ S$. This induces an involution
$``\;^*\;":C(A)\ra C(A)$ on $C(A)$. Let also $m_{ _A}(\ch,\;\mu)$ be the usual multiplicity form on
$C(A)$. Recall that  if $M$ and $N$ are $A$-representations affording characters $\ch, \mu \in C(A)$ respectively, then $m_{ _A}(\ch,\;\mu)$ is defined by $m_{ _A}(\ch,\;\mu):=\dim_{\kk}\Hom_{A}(M, N)$. We will use the notation $G(A)$ for the set of grouplike elements of $A$. \mdn For any right $A$-comodule $M$
with comodule structure $\rho:M \ra A\ot M$ denote by $C_M$ the subcoalgebra of
coefficients of $M$. Recall that $C_{M}$ is the smallest subcoalgebra $C\subset A$ with the property
that $\rho(M)\subset C\ot M$, see \cite{Lar}. If $d \in C(A^*)$ is the character of $M$ as an
$A$-comodule then $C_M$ is also denoted by $C_d$. 
\mdn
A {\it left coideal subalgebra }of $A$ is a left coideal $L$ of $A$ (i.e a vector subspace with $\Delta(L) \subset A \ot L$) which is also a unitary subalgebra of $A$. A left coideal subalgebra $L$ of $A$ is called {\it normal} if it is closed under the left adjoint action $ad_l$ of $A$ on itself. Recall that $ad_l(x)(a)=x_1aS(x_2)$ for all $a, x \in A$. If $L$ is a left normal coideal subalgebra of $A$ then $AL^+=L^+A$ is a Hopf ideal and let $A//L:=A/AL^+$ be the Hopf algebra quotient.  Let also $\pi_L:A \ra A//L$ be the canonical Hopf projection.\mdn It is known that any fusion subcategory of $A$ is of the form $\rep(A//L)$ for some normal left coideal subalgebra $L$ of $A$ (see for instance \cite{nrm}.)

\subsection{Hopf kernels of representations of semisimple Hopf algebras}
Let $A$ be any semisimple Hopf algebra and $M$ a finite dimensional left $A$-representation affording a charcater $\ch$. By \cite[Proposition 1.2]{Bker} one has that $|\chi(d)| \leq
\chi(1)\eps(d)$ for any $d \in \Irr(A^{*})$. Moreover, \cite[Remark 1.3]{Bker} implies that $\ch(d)=\eps(d)\ch(1)$ if and only if  the subcoalgebra $C_{d}$ acts trivially on $M$.
\mdn
Recall \cite{Bker} that the {\it Hopf kernel }$\hker_A(M)$ is defined as Hopf subalgebra of $A$ generated by the set of all characters $d \in \Irr(A^{*})$ such that $\ch_{M}(d)=\ch_{M}(1)\eps(d)$, where $\ch_{M}$ is the character associated to $M$. It follows that the Hopf kernel $\hker_A(M)$ of $M$ coincides to the largest Hopf subalgebra of $A$ that acts trivially on $M$. \mdn
Let $K$ be a normal Hopf subalgebra of $A$. Since $K=\oplus_{x \in \irr(K^{*})}C_{x}$ it follows from above that
\beq\label{reason}
\Rep(A//K)=\cap_{x \in \Irr(K^*)}\rep(\hker_{A^{*}}(x)^{*}). 
\eeq
\subsection{Brauer's theorem for left kernels of representations}

Let $M$ be an $A$-module and let $\mtr{LKer}_{ _A}(M)$ be the {\it left kernel of $M$}. Recall \cite{gmj}  that  $\mtr{LKer}_{ _A}(M)$ is defined by:
\begin{equation}\label{L}
\mtr{LKer}_{ _{A}}(M)=\{a \in A|\; a_1\ot a_2m=a\ot m,\;\;\text{for all}\;\; m\in M\}
\end{equation}
 Then by \cite{gmj} it follows that $\LKER_A(M)$ is the largest left coideal subalgebra of $A$ that acts trivially on $M$. It is also a normal left coideal subalgebra and obviously $\hker_A(M)\subset \lker_A(M)$. 
 \md
Next theorem generalizes a well known result of Brauer in the representation theory of finite groups. 
\bn{thm}\cite[Theorem 4.2.1]{gmj}.\label{charofim} Suppose that $M$ is a finite dimensional module over a semisimple Hopf algebra $A$. Then \beq
<M>=\mtr{Rep}(A//\mtr{LKer}_A(M))
\eeq
where $<M>$ is the fusion subcategory of $\mtr{Rep}(A)$ generated by $M$.
\end{thm}
%We also 
%\red{ multiplicity form.}
%It can be checked that $SL^+$ is a two sided ideal of $S$ and therefore $S//L$ is an algebra. %Moreover $S//L$ is a right $A$-comodule algebra via $\rho(\bar{s})=s_1\ot \bar{s_2}$. To see that $\rho$ is well defined note that $$\rho(s(l-\eps(l)1))=s_1l_1\ot \bar{s_2(l_2-\eps(l_2)1)}+s_1l\ot s_2-\eps(l)s_1\ot s_2$$
Suppose that $L$ is a normal left coideal subalgebra of a Hopf algebra $A$ and let $B:=A//L$ be the quotient Hopf algebra of $A$ via $\pi_L:A\ra B$. Then under the dual map $\pi_{L}^{*}$ it follows that $B^{*}$ is a Hopf subalgebra of $A^*$ that can be identified with:
\begin{equation}\label{forml}
\pi_{L}^*(B^*)=\{f\in A^*|f(al)=f(a)\eps(l)\;\; \text{for all}\;\; a\in A,\;l\in L\}.
\end{equation}

\bp\lb{hopfprod}
Let $B$ and $B'$ be two Hopf subalgebras of a semisimple Hopf algebra $A$. If $C(A^*)$ is commutative then
$$\dim_{\kk} (<B,B'>)=\frac{(\dim_{\kk} B)(\dim_{\kk} B')}{\dim_{\kk} (B\cap B')}.$$ %\blue{definition of lambda's}
\ep
\bpf Let $\Lam_{B}, \Lam_{B'}$ be the idempotent integrals of $B$ and $B'$ respectively.
Since $C(A^*)$ is commutative one has that $\Lam_B\Lam_{B'}=\Lam_{B'}\Lam_B$. This implies that $\Lam_B\Lam_{B'}=\Lam_{<B,\;B'>}$. The equality now follows counting the multiplicity of $1$ in both sides of the previous equality.
\epf
\br
Alternatively the statement follows from  \cite[Lemma 3.5]{DGNOI} applied to $\mathcal{C}=\rep(A^{*})$.
\er
If $S$ and $R$ are two left coideal subalgebras of $A$ denote by $SR$ the vector subspace of $A$ generated by elements of the type $sr$ with $s \in S$ and $r\in R$. Clearly $SR$ is a coideal of $A$ but not a coideal subalgebra unless $SR=RS$. Denote by $<S, R>$ the left coideal subalgebra of $A$ generated by $S$ and $R$. 

\bt\label{interms}
Let $L$ and $K$ be two normal coideal subalgebras of $A$. Then the following equalities hold in $A^*$:
\bn{enumerate}
\item\lb{int}
$(A//L)^*\cap (A//K)^*=(A//LK)^*.$

\item
$<(A//L)^*,\;(A//K)^*>=(A//(L\cap K))^*.$%\blue{this is used in the proof. can I shorten here, just the statement at the categorical level and then proof for the Hopf subalgebra generated.}
\end{enumerate}
\md
If the Grothendieck ring $\G_0(A)$ is commutative then
\beq
\dim_{\kk} (LK)=\frac{ \dim_{\kk} (L\cap K)}{(\dim_{\kk} L)(\dim_{\kk} K)}.
\eeq
\et

\bpf 
1.) It follows directly from Equation \eqref{forml}.
\md
2) Using Equation \eqref{forml} it is easy to see that 
\beqn<(A//L)^*,\;(A//K)^*>\subseteq (A//(L\cap K))^*.\eeqn

On the other hand, from  Takeuchi's correspondence \cite{tak} between Hopf subalgebras of $A^{*}$ and coideal subalgebras of $A$ it follows that there is a normal left coideal subalgebra $M(L,K)$ of $A$ such that
\begin{equation}\label{fl}
<(A//L)^*,\;(A//K)^*>=(A//(M(L,\;K))^*.
\end{equation}
Since $(A//L)^*\subseteq (A//M(L,\;K))^*$ it follows that $M(L,\;K)\subseteq L$. Similarly $M(L,\;K)\subseteq K$ and therefore $M(L,\;K)\subseteq L\cap K$. This implies that $\dim_{\kk} (<(A//L)^*,\;(A//K)^*>) \geq \dim_{\kk}(A//(L\cap K))^*)$. Then Equation \eqref{fl} implies the equality
\beqn<(A//L)^*,\;(A//K)^*>=(A//(L\cap K))^*.\eeqn
\md Suppose now that $\G_0(A)$ is commutative.
Using the previous proposition it follows that \beqn \dim_{\kk} <(A//L)^*,\;(A//K)^*> =\frac{\dim_{\kk} (A//K)^*\dim_{\kk} (A//L)^*}{\dim_{\kk} ((A//L\cap K)^*)}.\eeqn Since $<(A//L)^*,\;(A//K)^*>=(A//L\cap K)^*$ the equality follows.
\epf
For a normal left coideal subalgebra $L$ of $A$ it follows  that $M \in \Rep(A//L)$ if and only if $L$ acts trivially on $M$. Recall that $L$ acts trivially on $M$ if and only if $xm=\eps(x)m$ for all $x \in L$ and all $m \in M$.%Previous theorem can be restated as following:
\bc\label{lattic}
Let $L$ and $K$ be two normal left coideal subalgebras of a semisimple Hopf algebra $A$. Then one has the following:
\bn{enumerate}
\item
$\Rep(A//L)\vee \Rep(A//K)=\Rep(A//L\cap K).$
\item
$\Rep(A//L)\cap \Rep(A//K)=\Rep(A//LK).$
\end{enumerate}
\ec

\section{Kernels and centralizers in factorizable Hopf algebras}\lb{kcc}
%\subsection*{Quasitriangular Hopf algebras}
Recall that a Hopf algebra $A$ is called {\it quasitriangular} if $A$ admits an $R$-matrix, i.e. an element $R \in A\ot A$ satisfying the following properties:\\
$1) \;
R \Delta(x)=\Delta^{\cop}(x)R
$ for all $x \in A$.
\\$2)\;
(\Delta \otimes \id)(R)= R^{1} \ot r^{1} \ot R^{2}r^{2}
$
\\$3)\;
( \id \otimes \Delta)(R)=R^{1}r^{1} \ot r^{2}\ot R^{2}.
$
\\$4)\;
( \id \otimes \eps)(R)=1=(\eps \ot \id)(R).
$

Here $R=r=R^{1}\ot R^{2}=r^{1}\ot r^{2}$.
\md
%Then $u:=S(R^{2})R^{1}$ is the Drinfeld element of  $A$ and $S^2(a)=uau^{-1}$ for all $a\in A$. Since $A$ is semisimple one has $S^{2}=\id$ and therefore $u$ is a central element of $A$. Moreover $u^{-1}=\sum_{i=1}^mS^{-2}(R^{2})R^{1}$ and $\Delta(u)=(u \ot u)(R_{21}R)$, see \cite{dr}.
If $(A, R)$ is a quasitriangular Hopf algebra then the category of representations is a braided fusion category with the braiding given by
\beq
c_{M, N}:M\ot N\ra N\ot M,\; m\ot n \mapsto R_{21}(n\ot m)
\eeq
for all $M, N\in \co(\rep(A))$ (see \cite{Kas}). Recall that $R_{21}:= R^{2}\ot R^{1}$.
%\subsection{Factorizable Hopf algebras}
\md
A quasitriangular Hopf algebra $(A,\;R)$ is called {\it factorizable} if and only if the linear map  
\beq\label{phir}
\phi_A:A^* \ra A, \;\;f \mapsto (\id \ot f)(R_{21}R)\eeq is an isomorphism of vector spaces. In this situation, following \cite[Lemma 2.2]{schfact} $\phi_A$ maps the character ring  $C(A)$ onto the center of $Z(A)$ of $A$. Moreover by \cite[Theorem 2.1]{schfact} one  has that
\beq\lb{eqs2}
 \phi_A(f\chi)=\phi_A(f)\phi_A(\chi)
\eeq
for all $f \in A^*$ and $\chi \in C(A)$.  Thus $\phi_{A}|_{C(A)}:C(A)\ra Z(A)$ is an isomorphism of $\kk$-algebras.%Note also that by \cite[Theorem 2.3]{schfact} the map  $\phi_{R}$ restricted to $G(A)^{*}$ induces a group isomorphism $G(A^*) \ra G(A) \cap Z(A)$.

 \subsection{On the $S$-matrix for a factorizable Hopf algebra}
Let $A$ be a semisimple factorizable Hopf algebra. By \cite{EG} (see also \cite{KO02}) one has that the $S$-matrix $(s_{ij} )$ of the modular tensor category $\Rep(A)$ is given by \beq\label{s}s_{ij}=\chi_{i}(\phi_{A}(\ch_{j^{*}}))\eeq Moreover it is not difficult to see that for all $1\leq i,j \leq s$, one has that $s_{ij} = s_{ji}$, and $s_{ij} = s_{i^*j^*}$ and
$s_{ij^*} = s_{ji^*}$ (cf. \cite{BaKi, schfact}). Note that in this case inequality \eqref{is} becomes \beq\label{ifs}|s_{ij}|\leq \ch_{i}(1)\ch_{j}(1).\eeq

\subsection{Central primitive idempotents of the character ring}
Let $A$ be a semisimple Hopf algebra and $ \ch_{1}, _{\cdots}, \ch_{s}$ be the irreducible characters of $A$. Let also $e_{1}, _{\cdots},  e_{s}$ be their associated central primitive  central idempotents in $A$ and let $E_{j}:=\phi_{A}^{-1}(e_{j})$ for all $1 \leq j \leq s$. Since the map $\phi_A$ from Equation \eqref{phir} is an algebra isomorphism between $C(A)$ and  $Z(A)$ it follows that $\{E_j \}_{j=1}^{s}$ is a set of primitive central idempotents in the character ring  $C(A)$ of $A$. 

\bp\label{ineq} Let $A$ be a semisimple Hopf algebra and $W$ be an irreducible right $A$-comodule with associated character $d \in \irr(A^*)$.  Then for any primitive central idempotent $E_{j}\in C(A)$ one has that $E_j(d) \in \Z_{\geq 0}$ and $\sum_{j=1}^{s}E_{j}(d)=\eps(d).$ \ep 

  \bpf Note that $W$ can also be regarded as a left $A^{*}$-module. Its character, $\hat d:A^{*} \ra k$, is given by evaluation at $d$, i.e., $\hat d(f)=f(d)$ for all $f \in A^{*}$.  Since $C(A)$ is a semisimple subalgebra of $A^{*}$ it follows that $W$ can also be regarded as a left $C(A)$-module. 
\mdn 
Let $M_{j}$ be the associated simple $C(A)$-module corresponding to the central primitive idempotent $E_{j}\in C(A)$. Then \beqn E_{j}(d)=\hat d (E_{j})=m_{C(A)}(M_{j},\;W\dw^{A^{*}}_{C(A)}) \dim_{\kk} M_{j} \in \Z_{\geq 0}\eeqn where $m_{C(A)}(M_{j},\;W\dw^{A^{*}}_{C(A)})$ coincides to the multiplicity of the simple left $C(A)$-module $M_{j}$ associated to $E_{j}$ inside $W\dw^{A^{*}}_{C(A)}$. Clearly $\sum_{j=1}^{s}E_{j}(d)=\sum_{j=1}^{s}m_{C(A)}(M_{j},\;W\dw^{A^{*}}_{C(A)}) \dim_{\kk} M_{j}=\dim_{\kk} W=\eps(d).$ 
 \epf
For a subset $X \in \Irr(A)$ we denote by $<\ch\;|\;\ch \in X>$ the fusion subcategory of $\rep(A)$ generated by the irreducible modules $M$ whose characters satisfy $\ch_{M}\in X$. \br Note that for any Hopf subalgebra $L \subset A$ one has $\rep(L^{*}) \subset \rep(A^{*})$ via the canonical projection $A^{* }\ra L^{*}$.\er% or say you identify the Hopf kernel with a fusion subcategory of the dual. 

\bn{thm}\label{wding}
Let $A$ be a semisimple factorazible Hopf algebra and $d \in \Irr(A^*)$. Then one has the following equality:
\begin{equation*}
   \rep(\mtr{HKer}_{A^*}(d)^{*})'=<\ch_{j}\;|\;E_j(d)\neq 0>.
\end{equation*}
%where $j \mapsto \phi(j)$ is the above defined bijection. %between the set of primitive idempotents of $C(A)$ and the set of irreducible characters of $A$. 
\end{thm}
\bn{proof} Using the above remark, since $\mtr{HKer}_{A^*}(d)$ is a Hopf subalgebra of $A^{*}$ it follows that $\rep(\mtr{HKer}_{A^*}(d)^{*})$ is a fusion subcategory of $\rep(A)$. Moreover, by its definition one has that $\ch \in \rep(\mtr{HKer}_{A^*}(d)^{*})$ if and only if $\ch(d)=\ch(1)\eps(d)$.
\md
Using Equation \eqref{s} and \cite[Remark 3.4]{schfact} one can write that
\begin{equation}\label{ch}
 \ch_i=\sum_{j=1}^s\frac{s_{i{j}}}{s_{{0j}}}E_j.
\end{equation}
for all $1\leq i \leq s$.
On the other hand Proposition \ref{ineq} and inequality \eqref{ifs} gives that
\begin{equation*}
    |\ch_i(d)|=|\sum_{j=1}^s\frac{s_{ij}}{s_{0j}}E_j(d)|\leq \sum_{j=1}^s|\frac{s_{ij}}{s_{0j}}|E_j(d)\leq \ch_i(1)\eps(d)
\end{equation*}
for any character $d \in \Irr(A^{*})$.
Therefore one has that $\ch_i(d)=\eps(d)\ch_i(1)$ if and only if $s_{ij}=\ch_i(1)\ch_{j}(1)$ for all $j$ with $E_j(d)\neq 0$. It follows that $\ch_i \in \HKer_{A^*}(d)$ if and only if $\ch_i$ centralize any $\ch_{j}$ with the property that $E_j(d)\neq 0$. This completes the proof of the theorem.
\end{proof}
%\blue{Note that the above formula is at the level of characters.} %\blue{For a based subring $B\subset C(A)$ we denoted by $B'$ the based subring generated by all the irreducible characters that centralize all the irreducible characters contained in $B$.}
\bn{cor}\label{norm'}
Let $A$ be a factorizable Hopf algebra and $K$ a normal Hopf subalgebra $K$ of $A$. If $\Lam_{K}\in K$ is the idempotent intergral of $K$ then
\beq
\Rep(A//K)'=<\ch_{j}\;|\;E_j(\Lam_{K})\neq 0>.
\eeq
\end{cor}
\bpf
By Equation \eqref{reason} one has
$\Rep(A//K)=\cap_{x \in \Irr(K^*)}\rep(\hker_{A^{*}}(x)^{*})$. Therefore \beqn
\Rep(A//K)'=\vee_{x \in \Irr(K^*)}\rep(\hker_{A^{*}}(x)^{*})'=\vee_{x \in \Irr(K^*)}<\ch_{j}\;|\;E_j(x)\neq 0>.
\eeqn
\epf
By \cite[Proposition 4.1]{Lar}) note that  the idempotent integral of $K$ satisfies  \beq \Lam_{K}=\frac{\sum_{x \in \irr(K^{*})}\eps(x)x}{\dim K}\eeq  

Then Proposition \ref{ineq} shows that $E_{j}(\Lam_{K})\neq 0$ if and only if $E_{j}(x)\neq 0$ for some $x \in \Irr(K^{*})
$. Thus \beq
\Rep(A//K)'=<\ch_{j}\;|\;E_j(\Lam_{K})\neq 0>.
\eeq
\subsection{Kernels of grouplike elements}
Let $A$ be a semisimple factorizable Hopf algebra and $g \in G(A)$ be a grouplike element. Then Theorem \ref{centraliz} implies that 
\beq
\rep(\hker_{A^*}(g)^{*})'=<\ch_{g}>,
\eeq
where the character $\ch_{g}\in \Irr(A)$ is defined as follows.
First note that by Proposition \ref{ineq} there is a unique central primitive idempotent $E_g \in C(A)$ that does not vanish on $g$.  Then the corresponding irreducible character of the primitive idempotent $\phi_{A}(E_{g})$ is denoted by $\ch_{g}$.
\mdn Note that by \cite[Equation (1.3)]{Bd} the idempotent $\phi_{A}(E_{g})$ is a scalar multiple of 
\beq
F_g=\sum_{i=1}^{s} \ch_{i}^*(g) \ch_{i}. 
\eeq

\section{Kernels and centralizers in semisimple quantum doubles}\lb{kccd} In this section we will describe the M\"uger centralizer of some fusion subcategories of $\Rep( D(A))$ for any semisimple Hopf algebra $A$. In particular we will prove the two main results mentioned in the introduction. \md
Recall that the Drinfeld double $D(A)$ of a Hopf algebra $A$ is defined
by $D(A) \cong A^{*cop} \otimes A$ as coalgebras with the multiplication given by \beq (g
\bowtie a)(f \bowtie b)=\sum g(a_1\rightharpoonup f \leftharpoonup S^{-1}a_3)\bowtie
a_2b,\eeq for all $a,b \in A$ and $f,g \in A^*$. Moreover its antipode is given by $S(f
\bowtie h)=S^{-1}(h)S(f)$. Recall that $(a\lh f \rh b)(x)=f(bxa)$ for all $a,b \in A$ and $f \in A^{*}$. It is well known that $D(A)$ is a semisimple Hopf algebra if
and only if $A$ is a semisimple Hopf algebra \cite{Montg}. % where $\mtr{Rep}(A)^{\mtr{rev}}:=\rep(A^{\mtr{op}})$.
\md
If $A$ is a semisimple Hopf algebra and $\cC=\Rep(A)$
 then  $\mathcal{Z}(\cC)\simeq \Rep(D(A))$  where $D(A)$ is its quantum double \cite{Kas}. 
\md
Moreover, the Drinfeld double $D(A)$ is a quasitriangular Hopf algebra with $R$-matrix: $R=\sum_{i=1}^n(\epsilon \bowtie b_i) \ot( b_i^* \bowtie 1)$ where $\{b_{i}\}$ is a vector basis for $A$  and $\{b_{i^{*}}\}$. %Moreover in this case $u=\sum_{i=1}^nS^{-1}b_i^* \bowtie b_i$ and $u_D^{-1}=\sum_{i=1}^nS^{2}b_i^* \bowtie b_i$. 
In this case one has 
 $R_{21}R=\sum_{i,j=1}^n( b_j^* \bowtie b_i)\otimes (\epsilon \bowtie b_j)(b_i^* \bowtie 1) $ and by \cite[Lemma 1.1]{EG} $D(A)$ is a factorizable Hopf algebra. The linear map $\phi_{D(A)}:D(A)^{*}\ra D(A)$ is given by \beqn f \mapsto \sum_{i,j=1}^n( b_j^* \bowtie b_i)f[(\epsilon \bowtie b_j)(b_i^* \bowtie 1) ].\eeqn
%\blue{Note that under the canonical identification of vector spaces $D(A)^{*}\cong A\ot A^{*}$ the map $\phi_{R}$ satisfies $\phi_{R}(h\ot f)=(\epsilon \ot h)(f\ot 1)=hf$. }Also,
%if $f=\sum_j h_j\ot f_j \in C(D(H))$ then $\phi_{R}(f)=\sum_jh_j\bowtie f_j$ (see \cite{KSZ}) \bb{DO I use it?}.
%}For any vector subspace $S$ of $D(A)$ denote by $\cD(S)$ the subcategory of $\cD$ consisting of those objects on which $S$ acts trivially.

\subsection{Left Kernels of normal Hopf subalgebras}\label{mod}
 Recall \cite{Zind} that $A$ can be regarded as a $D(A)$-module via the action
\beq\label{ms}
(f\bwt a).b=(a_{1}bS(a_{2}))\lh Sf
\eeq
for all $a, b \in A$ and all $f \in A^{*}$.
Then note that any normal Hopf subalgebra $K$ of $A$ can be regarded as $D(A)$-submodule of $A$.  %We will use Brauer's theorem \ref{charofim} to find the fusion subcategory $<K>$ of $\mtr{Rep}(D(A))$ generated by the $D(A)$-submodule $K\subseteq A$.

\bp\lb{leftkr} Let $K$ be a normal Hopf subalgebra of $A$. Then
$$\LKER_{\Dr(A)}(K) \supseteq (A//K)^*\bowtie \LKER_A(K),$$
where $\LKER_A(K)$ is the left kernel of $K$ regraded as a $A$-module via the adjoint action.
\ep
\bpf
From the module structure of Equation \eqref{ms} it can easily be seen that the left coideal subalgebra  $(A//K)^*\bowtie \LKER_A(K)$ of $D(A)$ acts trivially on $K$. Indeed $(A//K)^{*}$ acts trivially on $K$ due to Equation \eqref{forml}. Since $\LKER_{\Dr(A)}(K)$ is the largest coideal subalgebra of $D(A)$ that acts trivially on the $D(A)$-module $K$ the inclusion $(A//K)^*\bowtie \LKER_A(K)\subseteq \LKER_{\Dr(A)}(K)$ follows.
\epf
%\br Let $K$ be a normal Hopf subalgebra of $A$. Regard $K$ as an $A$-module via the adjoint action. Then by Remark \ref{hadj} it follows that $\HKer_A(K)$ is the largest Hopf subalgebra of $A$ that commutes with $K$. 
 %\er
\subsection{Primitive central idempotents in $C(A)$}
There is another realization (see \cite[Proposition 6.3]{KSZ}) of $A$ as a $D(A)$-module on the vector space $A^{*}$. This realization coincides with the trivial $A$-module induced to $D(A)$ and has the following structure:
\beqn
(f \bwt a)g=f(a_{1}\rh g \lh S(a_{2}))
\eeqn
for all $f, g \in A^{*}$ and $a \in A$. Then by  \cite{KSZ} each homogenous $D(A)$-component of $A^{*}$ can be described as $A^{*}E_{i}$ for some central primitive idempotent $E_{i }\in C(A)$, see \cite{KSZ}. Moreover $\{E_{i}\}_{i=1,s}$ is a complete set of central orthogonal primitive idempotents of $C(A)$.
\mdn
 Recall by \cite[Proposition 3.1]{repalg} that the Fourier transform $\mathcal{F}:A \ra A^{*}$ given by $a \mapsto a \rh t$ is a morphism of $D(A)$-modules where $t \in A^{*}$ is the idempotent integral of $A^{*}$.
\mdn
Let $A=V_{1}\oplus V_{2} \oplus_{\cdots} \oplus V_{s}$ be the decomposition of $A$ into homogenous $D(A)$-components and suppose that $\mathcal{F}(V_{i})=A^{*}E_{i}$. Define the linear functionals $p_i \in A^*$ by the relation that $p_i$ coincides to $\eps_A$ on $V_i$ and vanishes on any other $V_j$ with $i\neq j$. Then by \cite[Theorem 5.13]{repalg} it follows that $p_{V_{i}}=S(E_{i})$ if $\mathcal{F}(V_{i})=A^{*}E_{i}$. 
\bl
Suppose that $K$ is a normal Hopf subalgebra of a semisimple Hopf algebra $A$ and $d \in \Irr(K^{*})$ is a cocommutative element of $K$. With the above notations if $p_{V_{i}}(d)\neq 0$ then $V_{i} \subset K$.
\el
\bpf
If $K$ is a normal Hopf subalgebra of $A$ then $K$ is a full isotopic submodule of $A$, see \cite[Proposition 5.4]{repalg}.\
Thus  without loss of generality one may suppose that as homogenous $D(A)$-components one has the decomposition $ K=V_{1} \oplus {\dots} \oplus V_{r}\; \text{and}\; A=K \oplus V_{r+1} \oplus {\dots} \oplus V_{s}$ for some $1\leq r \leq s$.
Then clearly if $p_{V_{i}}(d)\neq 0$ then $ 1\leq i \leq r$ and therefore $V_{i} \subset K$. 
\epf

\subsection{Proof of Theorem \ref{justonecomp}}
Now we are ready to prove Theorem \ref{justonecomp}
% Let $V$ be a simple $D(A)$-submodule of $A$. Then by KSZ $E_{V}\dw^{D(A)}_{A}=q_{W}$ where $W$ is the homogenous component containing $V$.
%
%\bb{do I need commutative ring here? Is this the bijection described by Sommerhauesr?}
\bpf
Let $\ch_{1},{\dots}, \ch_{l}$ be the irreducible characters of $D(A)$ and $\tilde{e}_{1}, {\dots},\tilde{e}_{l}\in \cZ(D(A))$ be their associated central idempotents in $D(A)$. Denote also by $\tilde{E}_{j}=\phi_{D(A)}^{-1}(\tilde{e}_{j})$ the primitive central idempotents of $D(A)$.
\md
Similar to Equation \eqref{reason} one has that $\cD(K)=\cap_{d \in \Irr(K^*)}\rep(\HKer_{D(A)^*}(d)^{*})$. Thus 
\beanon
\cD(K)' & = &  (\cap_{d \in \Irr(K^*)}\rep(\HKer_{D(A)^*}(d))^{*})'=\vee_{d \in \Irr(K^*)}\rep(\HKer_{D(A)^*}(d)^{*})'\eeanon
Since $D(A)$ is a factorizable Hopf algebra, Theorem \ref{wding} implies that
\beanon
\cD(K)'  &= & \vee_{d \in \Irr(K^*)}(<\ch_j\:|\;\tilde{E}_j(d)\neq 0>)\\&=&<\ch_{j}\;|\;E_{j}(d)\neq 0\text{\;for some\;} d \in\Irr(K^{*})>.
\eeanon
\md
Let $\tilde{E}_j$ be a primitive central idempotent of $C(D(A))$ such that $\tilde{E}_j(d)\neq 0$ for some $d \in \Irr(K^*)$. It follows that $\tilde{E}_j \in D(A)^{*}$ has a nonzero restriction to $A$. Therefore by \cite[Theorem 6.3]{KSZ} $\tilde{E}_{j}\dw^{D(A)}_{A}=E_{j}$ where $E_{j}$ is a central idempotent of $C(A)$ and $\ch_{j}$ is the character of any simple $D(A)$-submodule of the homogenous component $A^{*}E_{j}$. It follows by the above arguments that $S(E_{j})=p_{V_{j}}$ . Then the previous lemma implies that $\ch_j$ is the character of a simple $D(A)$-submodule $V_{j}$ of $A$ with $V_{j}\subset K$. Thus: \beanon
\cD(K)' & = & <\ch_{j}\;|\;E_{j}(d)\neq 0\text{\;for some\;} d \in \Irr(K^{*})>\\ & = & <V_{j}\;|\;V_{j}\subset K>=<K>.
\eeanon
\epf
It is well known that $D(A) \cong D(A^{*\;op\;cop})^{op}$ as Hopf algebras via $f \bwt a \mapsto a \bwt f$ (see for instance \cite[Theorem 3]{mintr}).
By duality this implies the following result:
\bn{cor} 
If $L$ is a normal Hopf subalgebra of $A$ then
\beqn \cD((A//L)^{*})'=<(A//L)^*>,\eeqn
where $<(A//L)^{*}>$ denotes the fusion subcategory of $\rep(D(A))$ generated by $(A//L)^{*}$, regarded as a $D(A)$-module via the above isomorphism $D(A) \cong D(A^{*\;op\;cop})^{op}$.
\ec
Next proposition is the dual version of the Proposition \ref{leftkr}.
\bp\label{leftkrdu}
Let $L$ be a normal Hopf subalgebra of $A$. Then
\beqn
\LKER_{D(A)}((A//L)^{*\cop})\supseteq \LKER_{A^{*}}((A//L)^{*})\bwt L
\eeqn
\ep
%\blue{Check if that also holds for normal coideal subalgebras. shouldn't be right ?? since the coposite for comutiplication.}
\newcommand{\sms}{\;semisimple\;}
\subsection{Proof of Theorem \ref{centraliz}}
\bpf
Suppose that $K$ and $L$ are normal Hopf subalgebras of $A$. Using Theorem \ref{justonecomp} and the previous corollary it follows that:
\beanon
\cD(K,\;L)' & = &  (\cD((A//K)^*)\bigcap \cD(L))'=\\ &=& \cD((A//K)^*)'\bigvee \cD(L)' =<(A//K)^*>\bigvee <L>
\eeanon
 Theorem \ref{charofim} implies that \beq
\cD(K,\;L)'= \rep(D(A)// \LKER_{D(A)}(A//K)^* ) \bigvee  \rep (D(A)// \LKER_{D(A)}(L )) 
 \eeq 
Then Corollary \ref{lattic} gives that 
\beq
\cD(K,\;L)'= \rep(D(A)//( \LKER_{D(A)}(A//K)^* \bigcap  \LKER_{D(A)}L)))
\eeq
Proposition \ref{leftkr} and its dual version, Proposition \ref{leftkrdu}, imply that 
\begin{eqnarray*}
&& \LKER_{D(A)}((A//K)^*) \bigcap  \LKER_{D(A)}(L) \supseteq  \\& \supseteq& ((A//K)^*\bwt \LKER_{A}(K))\bigcap ((A//L)^*\bwt \LKER_{A^{*}}((A//L)^{*})) 
\end{eqnarray*}
Note that since $B(K,L)=(A//K)^{*}\bwt L$ is a normal Hopf subalgebra of $D(A)$ it follows by \cite[Theorem 4.5]{nrm} that the pairs $(L,K)$ and $((A//L)^{*},\;(A//K)^{*})$ are commuting pairs of Hopf subalgebras of $A$ and $A^{*}$ respectively. Since $L$ commutes elementwise with $K$ it follows that $L$ acts trivially on the $A$-representation $K$ and therefore $\LKER_{A}(K)\supseteq L$. Similarly $\LKER_{A^{*}}((A//L)^{*})\supseteq (A//K)^{*}$. These two inclusions imply that:
\begin{eqnarray*}
 ((A//K)^*\bwt \LKER_{A}(K))\cap ((A//L)^*\bwt \LKER_{A^{*}}(A//L)^{*})  =(A//K)^{*} \bwt L
\end{eqnarray*}
and $\rep(D(A)//( \LKER_{D(A)}(A//K)^* \bigcap  \LKER_{D(A)}L)))\subseteq \cd(K,L)$.
Hence one can conclude that $\cD(K,\;L)' \subseteq \cD(L,\;K)$. \md On the other hand, note that \beqn\fp \;\cd(K, L)=\frac{(\dim_{\kk} A )(\dim_{\kk} L)}{\dim_{\kk} K}.\eeqn Thus \beqn \fp \;(\cd(K, L)')=\frac{(\dim_{\kk} A) ^{2}}{\fp\; \cd(K, L)}=\frac{(\dim_{\kk} A )(\dim_{\kk} K)}{\dim_{\kk} L}=\fp \;\cd(L, K).\eeqn
This implies that  $\cd(K, L)'=\cd(L, K)$.% which implies the con}
%\blue{need to verify that these are normal coideal subalgebras. }
\epf

\bc If $K$ is a normal commutative Hopf subalgebra of a semisimple Hopf algebra $A$ then $\cd(K)'\subseteq\cd(K)$.\ec

\bpf
Note that $K$ regarded as $D(A)$-representation via Equation \eqref{ms} satisfies $K\in \cd(K)$ since $K$ is a commutative Hopf algebra.  Therefore Theorem \ref{centraliz} implies that $\cd(K)'=<K>\subseteq \cd(K)$.%\bb{comodule structure is a problem}
\epf
\subsection{Normal Lagrangian subcategories of the category of representation of a Drinfeld double}

We call a fusion subcategory of $\rep(A)$ normal if it is of the type $\rep(A//L)$ for a normal Hopf subalgebra  $L$ of $A$. This definition agrees with the definition of a normal fusion subcategory from \cite{tensor-exact}. Let $\cc$ be a premodular category with braiding $c$ and twist $\teta$. According to \cite{DGNO}, a
fusion subcategory $\ce$ of $\cc$ is called {\it isotropic } if $\teta$ restricts to the identity on $\ce$, i.e.,
if $\teta(X)=\id_{X}$ for all $X\in \ce$. Moreover, an isotropic subcategory $\ce$ is called {\it Lagrangian} if $\ce'=\ce$.
Recall also that a fusion category $\cc$ is called {\it hyperbolic}  if it has a Lagrangian subcategory and $\cc$ is called
{\it anisotropic} category if it has no non-trivial isotropic subcategories.
%andanisotropic if it has no non-trivial isotropic subcategories.}
\bt
Suppose that $F,G$ are finite groups and $A$ is an abelian extension fitting the extension
\beq
\kk \ra \kk^{G } \ra A \xra{\pi_{1}} \kk F\ra \kk
\eeq
Then $\cd(\kk^{G}, \kk^{G})$ is a Lagrangian subcategory of $D(A)$.
\et
\bpf Note that $A^{*}$ fits the exact sequence
\beq
\kk \ra \kk^{F } \ra A^{*} \xra{\pi_{2}} \kk G\ra \kk.
\eeq Since $A/\kk^{G}\cong \kk^{F}$ it follows that $B(\kk^{G}, \kk^{G})$ is a normal Hopf subalgebra of $D(A)$. Moreover, 
by Theorem \ref{centraliz} one has that $\cd(\kk^{G}, \kk^{G})'=\cd(\kk^{G}, \kk^{G})$. Let $\pi:D(A)\ra D(A)/B(\kk^{G}, \kk^{G})$ be the canonical Hopf projection.
\mdn  Using \cite[Poposition XIV.2]{Kas} note that the twist structure on $D(A)$ is given by the left multiplication by the Drinfeld element $u$ associated to $A$, see also \cite{EG}.
\mdn
On the other hand, in our case, the Drinfeld element of $D(A)$ can be written as
\beqn
u=\sum_{x\in F, a\in G}S(q_{x} \# a) (p_{a} \# x)
\eeqn
where $\{q_{x}\}$ and $\{p_{a}\}$ are the dual group element bases in $\kk^{F}$ respectively $\kk^{G}$.
It follows that $\pi(u)=\sum_{a,x}\pi_{2}(q_{x} \# a)\pi_{1}(p_{a} \# x)=1$, which shows that $\pi$ acts as identity on each object of $\cd(\kk^{G}, \kk^{G})$.
\epf
\br Note that \cite[Theorem 4.5]{DGNO} shows that in the situation of an abelian extension $A$ one has that
$\rep(D(A))$ is a hyperbolic modular category, i.e. it is braided tensor equivalent to the center $\cz(Vec^{\omega}_{G})$ for some finite group $G$ and some $\omega \in H^{3}(G, \kk^{*})$. In this way one can recover \cite[Theorem 1.3]{N1}.
\er
Recall the subcoalgebra of coefficients $C_M$ associated to a any right $A$-comodule $M$. By duality, any left $A$-module $V$ with associated character $\ch \in C(A)$ can be regarded as a right $A^*$-module and one can associate to it its subcoalgebra of coefficients $C_{\ch}\subset A^*$.
\subsection{Normal Hopf subalgebras of $D(A)$}
Let $L$ be a normal Hopf subalgebra of $A$. Recall that an irreducible character $\al$ of $L$ is called $A$-stable if there is a character $\ch \in \mtr{Rep}(A)$ such that 
$\ch\dw^{ ^A}_{ _L}=\frac{\ch(1)}{\al(1)}\al$. Such a character $\ch \in \Irr(A)$ is said to seat over the character $\al \in \Irr(L)$. The set of all irreducible $A$-characters seating over $\al$ is denoted by $\Irr(A|_{\al})$. Denote by
$G_{A}^{st}(L)$ the set of all $A$-stable linear characters of $L$. Clearly $G_{A}^{st}(L)$ is a subgroup of the group of grouplike elements $G(L^*)$ of the dual
Hopf algebra of $L$.
\mdn
Suppose that $K$ and $L$ are two normal Hopf subalgebras of $A$
and let $G$ be a finite group that can be simultaneously embedded in  $G^{st}_A(K)$ and
$ G^{st}_{A^*}((A//L)^*)$ via the emebddings $\psi_1:G \hookrightarrow G^{st}_A(K)$
and respectively $\psi_2:G \hookrightarrow G^{st}_{A^*}((A//L)^*)$. Let $B(K, L,G, \cX, \psi_{1}, \psi_{2})$
be the subcoalgebra of $D(A)$ defined by \beq B(K, L,G, \cX, \psi_{1}, \psi_{2})=\bigoplus_{x \in
G}C_{\psi_1(x)\uw^{ ^{A}}_{ _{K}}}\bowtie C_{ \psi_2(x)\uw^{ ^{A^*}}_{ _{(A//L)^*}}}.
\eeq 
\mdn Recall by \cite[Theorem 4.2]{nrm} that any normal Hopf subalgebra of $D(A)$ is of the type $B(K,L, \cX, \psi)$ where $K, L$ are normal Hopf subalgebras of $D(A)$ and $\cX$  is a finite group satisfying the above properties, for more details see \cite[Subsection 3.2 ]{nrm}. \mdn %Clearly one has that $B(K, L) \subseteq  B(K, L,G, \cX, \psi_{1}, \psi_{2})$.  
Let  $\cd(K, L, \cX,\psi)$ be the category of representations of the quotient Hopf algebra $D(A)//B(K,L,\cX,\psi)$.
\bp Let $A$ be a semisimple Hopf algebra.
Then $\rep(D(A))$ has a normal Lagrangian fusion subcategory if and only if $A$ is an abelian extension Hopf algebra.
\ep
\bpf
Suppose that $\cd(K,L, \cx, \psi)$ is a Lagrangian fusion subcategory of the category $\rep(D(A))$. In particular one has that $\cD(K, L, \cx, \phi)'=\cD(K,L, \cx, \psi)$. Note that the Hopf algebra inclusion $B(K,L)\subset B(K,L, \cx, \psi)$ induces an inclusion of fusion subcategories $\cD(K, L)\supseteq \cD(K,L, \cx, \psi)$. Thus $\cD(L,K)=\cD(K, L)'\subseteq D(K, L, \cx, \phi)'=D(K,L, \cx, \phi)$. This inclusion implies that $B(L, K)\supseteq B(K, L, \cx, \phi)$. Since $B(K, L)\subseteq B(K,L, \cx, \psi)$ this shows that $B(L, K)\supseteq B(K,L)$, and  therefore $K=L$. By \cite[Theorem 4.3]{nrm} one has that $[K,L]=0$, i.e. $K$ is a commutative Hopf algebra. Similarly $(A//K)^{*}$ is commutative and therefore $A$ is an abelian extension.
 \epf
It would be an interesting question to determine M\"uger's  centralizer for any normal fusion subcategory of $\rep(D(A))$ and to decide if the centralizer is also a normal fusion subcategory.
\bibliographystyle{amsplain}
\bibliography{hbc}
 \ed
  \br
 This proposition shows that in general $B(L, K)\supseteq B(K', L, \psi')$.
 \er
\section{Fusion subcategories of Drinfeld centers $\cZ(\cA)$ and their centralizers}\lb{qst}
\blue{Two questions on the formule kernels la for left  kernels and $\cd(L)'$}.\mdn
Let $\ca$ be a fusion category and let $\cZ(\cA)$ be \blu{its Drinfeld center} with forgetful functor $F : \cZ(\cA) \ra  \ca$. Let $\cc \subset \cz(\ca)$ be a fusion subcategory and let $\cc' \subset \cz(\ca)$ be its M$\ddot{u}$ger centralizer in $\cz(\ca)$.\mdn
\blue{Definition of $\cc \cap I(1)$ and $\cA(\cc \cap I(1))$.}
\mdn
Note that the right adjoint functor $I:\ca \ra Z(\ca)$ of $F$ defines a tensor equivalence $\ca \xra{I} Z(\ca)_{I(1)}$. Next Theorem is  \cite[Theorem 3.15]{dno}.
\stat 1\mdn
\blue{Their correspondence}
\stat 2\mdn
\blue{Their theorem}\mdn
\bt \cite{dno}\lb{dmno}
Let $\cc$ be a fusion category of $\cZ(\cA)$. Then $\ca(\cc \cap I(1))$ is precisely the image $F(\cc')$ of $\cc'$ in $\ca$ under the forgetful functor.
\et
\noindent
One has that $F(\cc)\cong \cc_{\cc\cap I(1)}$ since the adjoint functor of $F|_{\cc}:\cc\ra F(\cc)$ is given by $X \ra I(X)\cap \cc$. This implies that $\fp(F(\cc))=\frac{\fp(\cc)}{\fp(\cc \cap I(1))}$. \blue{this is explanation from the proof}
\mdn
\stat 3\mdn
Define $\cd(L)$ as those $D(A)$-modules that  are receiving trivial action from $L$.
\mdn
If $\ca = \rep(A)$ for a semisimple Hopf algebra $A$ and $\cc \subset \rep(D(A))\simeq \cz(\cA)$ then clearly $\cc \cap I(1)$ is the largest normal coideal subalgebra $L$ of $A$ such that $L \in \cc$. 
\bn{example} Suppose that $\ca=\rep(A)$ for a semisimple Hopf algebra $A$.
It is easy to see that if $\cc=F^{-1}(\rep(A//L))$ then $\cc'$ consists of those modules that are dyslectic with respect to $L$.\blue{this does not come from their correspondence }
\end{example}
\bp
With the above notations one has that $\cd(L)=F^{-1}(\Rep(A//L))$.
\ep
\bpf
Straightforward knowing that $\rep(A//L)$ consists of those $A$-modules receiving trivial action from $L$.
\epf
\stat 4\mdn
The previous theorem says that if $\cc\cap I(1)=L$ then all $D(A)$-modules which are objects of $\cc^{'}$ are trivial when restricted to $L$, thus $\cc'\subset \cd(L)$. This is equivalent to $\cc \supset \cd(L)'$.
\stat 5\mdn
\bab$F(\rep(D(A)//N))=\rep(A//N\cap A)$. \blue{is it true?}\eab
\stat 5'\mdn
\beqn
\fp(F(\cd(L)))=\frac{\fp(\cd(L))}{\dim M}
\eeqn
\stat 6 \mdn
\bab One has that $F(\cd(L))=\rep(A//(N(L)\cap A))$ where $N(L)$ is the smallest normal coideal subalgebra of $D(A)$ containing $L$.\eab
\stat 7\mdn
Suppose now that $\cc \cap I(1)\supseteq L$. Then by the previous remark one has that $\cc' \subseteq \cd(\cc \cap I(1))\subseteq \cd(L)$.
\stat 8\mdn\blue{Applying Theorem \ref{dmno} one has that $$F(\cd(L)')=\rep(A//M)=\;_{A}<L>.$$}
For $L=K$ a normal Hopf subalgebra of $A$ it follows that $$F(\;_{D(A)}<K>)=\;_{A}<K>.$$
\stat 9\mdn
\bp
One has that $\cd(L)\cap I(1)=\lker_{A}(L)$, the largest normal coideal subalgebra of $A$ that commutes elementwise with $S(L)$.
\ep
\bpf
One has that $\cd(L)\cap I(1)$ is the largest normal coideal subalgebra of $A$ that commutes with $S(L)$.
\mdn
Note that $M:=\lker_{A}(L)$ is  also the largest coideal subalgebra of $A$ that commutes with $S(L)$.
\epf\noindent
\stat 10\mdn
\bp Let $L$ be a normal left coideal subalgebra of $A$. 
With the above notations one has that 
\beq\label{incl}
\cd(L)'\subset \cd(M)
\eeq
where $M:=\lker_{D(A)}(L)\cap A$\ep
\bpf
Take $\cc=\cd(L)$. Then \blue{$\cd(L)\cap I(1) =M$.}\blue{Previous stament implies that $\cd(L)'\subseteq \cd(\lker_{A}(L))=\cd(M)$.}
\epf\blue{Compute $\fp$ dimensions. So far applied the weaker versions of the theorem.}

\bc
If $L$ commutes with $S(L)$ then $\cd(L)'\subseteq\cd(L)$.
\ec
\bpf
If $L$ commutes with $S(L)$ then clearly $L\subseteq M$ and therefore $\cd(L)'\subseteq\cd(M)\subseteq \cd(L)$. \epf
\stat 11
\bn{conj} \lb{coidconj}Let $L$ be a normal left coideal subalgebra of $A$. Then one has that \beq \lb{coidc}\cd(L)'=<L>\eeq the fusion subcategory generated by $L$.\end{conj}
\bpf\epf
\br
Clearly $<L>\subseteq \cd(M)$. 
\er
\stat 12\mdn
Note that for group algebras the above conjecture is verified.  \blue{\onh one does not have the equality $\cd(L)'=\cd(M)$.}

Indeed, $\cd(N)=\cs(C_{G}(N), N, 1)'=\cs(N, C_{G}(N), 1)$ while $\lker_{D(G)}(kN)\cap kG=kC_{G}(N)$ and $\cd(C_{G}(N))=\cs(C_{G}(C_{G}(N)), C_{G}(N), 1)$.
\blue{\bp
Let $L$ be a normal left coideal subalgebra of $A$. Then with the above notations one has that \beq \lb{scoidc}\cd(L)'\subset <L>\eeq
\ep
\bpf
Note that $\cd(L)\cap I(1)=\lker_A(L)$. Then the proof follows from Theorem \ref{dmno} and Brauer's Theorem since the pre image of $F$ by $\Rep(A//\lker_A(L))$ is $\Rep(D(A)//\lker_{D(A)}(L))$.
\epf
}
\bn{conj} $$F^{-1}(\rep(A//\lker_A(L))=<L>$$\end{conj}
\stat 13\mdn
Clearly $$F^{-1}(\rep(A//\lker_A(L)))=\cd(\lker_A(L))=\rep(D(A)//N(\lker_A(L))$$

\stat 15\mdn
\bn{conj} \lb{coidconj}Let $L$ be a normal left coideal subalgebra of $A$. Then one has that \beq \lb{coidc}\cd(L)'=<L>\eeq the fusion subcategory generated by $L$.\end{conj}
\blu{Conjecture: $\cd(L)'=\cd(\lker_A(L))\cap \cd((A//L)^*)$ verify it for Drinfeld doubles of Kac algebras.}
\stat 16 Group situation\mdn
\subsection{Irreducible reprsentations of $D(\kk G)$}

\subsection{Fusion subcategories of $\mtr{Rep}(D(G))$.} In this subsection we recall the parametrization of fusion subcategories of $\Rep(D(G))$ given in \cite{NNW}.
%_______________________________________1.8____________________________Let $\mtc{D}$ be a fusion subcategory of $\mtc{C}:=\mtr{Rep}(D(G))$. Then following
Cf. \cite{NNW} the fusion subcategory $\mtc{D}$ is completely determined by two canonical normal
subgroups $K_{ _{\cd}}$ and $H_{ _{\cd}}$ of ${G}$ and a
$G$-invariant bicharacter $B_{ _{\cd}}:K_{ _{\cd}}
\times H_{ _{\mtc{D}}} \ra \mathbb{C}^*$. The subgroups
$K_{ _{\mtc{D}}}$ and $H_{ _{\mtc{D}}}$ are defined as follows:
$K_{ _{\mtc{D}}} := \{gag^{-1}\; | g \in G\;\; \text{and} \;\;(a, \gamma)
\in \mtc{D} \;\;\text{for some} \;\;\gamma\}$ and $H_{ _{\mtc{D}}}$ is
the normal subgroup of $G$ such that $\mtc{D}\cap
\mtr{Rep}(G)=\mtr{Rep}(G/H_{ _{\mtc{D}}})$. Note that $K_{ _{\mtc{D}}}$ is
the fusion subcategory of $\mtr{Vec}^{G}$ determined by restricting all simple objects of $\mtc{D}$ to $\mathbb{C}G^*$.

The bicharacter %$B_{ _{\mtc{D}}}$
$B_{ _{\mtc{D}}} : K_{ _{\mtc{D}}}  \times H_{ _{\mtc{D}}} \ra \mathbb{C}^*$
is defined by
$B_{ _{\mtc{D}}}(g^{-1}ag, h) :=\frac{\gamma(ghg^{-1})}{\gamma(1)}$
if $(a, \gamma) \in \mtc{D}$. This is well defined and does not
depend on $\gamma$ by \cite{NNW}.
\md
Recall from \cite{NNW} that a such bicharacter is called $G$-invariant if and only if $B(xkx^{-1}, xhx^{-1})=B(k, h)$ for all $x \in G$, $k \in K$ and $h \in H$.
\mdn
Conversely, any two normal subgroups $K$ and $H$ of $G$ that
centralize each other elementwise together with a $G$-invariant bicharacter $B: K\times H \ra \mathbb{C}^*$ give rise to a fusion category denoted by $S(K, H, B)$ in \cite{NNW}. It is defined as the full
abelian subcategory of $\mtr{Rep}(D(G))$ generated by the objects
$(a ,\gamma)$ such that $a \in K \cap \mtc{R}$ and $\gamma \in
\mtr{Irr}(C_G(a))$ such that $\gamma(h)=B(a, h)\gamma(1)$ for all $h
\in H$.
Let $G$ be a finite group and $N$ be a normal subgroup of $G$. In the paper \cite{NNW} the authors have shown the following parametrization for fusion subcategories of $\Rep(D(G))$.
\blue{Make the conjecture as a question and write directly it is verified in the three cases, without mentioning details.}
\subsection{The conjecture for $A=\kk G$.}Next we will verify the previous conjecture for $A=\kk G$.
\bl For a normal subgroup $N$ of $G$ one has the following 
\bne \item \beqn
\lker_{\kk G}(\kk N)=\kk C_{G}(N),\eeqn where $C_{G}(N)$ is the centralizer of $N$.
\item \beqn
\cd(\kk N)=\cS(C_G(N), N, 1).\eeqn
\item \beqn\;_{D(\kk G)}<\kk N>=\cS(N, C_G(N), 1)\eeqn
\ene
\el
\bpf
\blue{Indeed..}
\epf
\noindent

\mdn
Then it follows by \cite[Theorem]{NNW} that $\cd(\kk N)'=\cS(N, C_G(N), 1)=<\kk N>$.
\md
\bn{example}
Note that if $N$ is a proper subgroup of the center $Z(G)$ then $C_G(C_G(N))=C_{G}(G)=Z(G) \neq N$. Thus in this situation one has that $\cd(\kk N)'\neq \cd(\kk C_{G}(N))$ showing thta the inclusion in Equation \eqref{incl} is strict.\end{example}
\md
\stat 17\mdn
Dual of group situation.\mdn
\subsection{On the Drinfeld double of  $\kk^{G}$}
Any left coideal subalgebra of $\kk^{G}$ is of the type
\beqn
\kk^{(G/M)_{l}}=\{f \in \kk^{G}\;|\; f(gm)=f(g)\;\text{for all\;}g \in G, m \in M\}
\eeqn
for some subgroup $M$ of $G$.
\mdn \blue{note that $D(\kk G^{\op})^{\op}\simeq D(\kk^{G})$.}
\bp
If $\kk^{(G/M)_l}$ is a left coideal subalgebra of $\kk^{G}$ then
\bne \item\beqn
\cd(\kk^{(G/M)_l})=\cs(\mtr{core}_G(M), \{1\}, 1).
\eeqn
\item\beqn
<\kk^{(G/M)_l}>=\cs( \{1\}, \mtr{core}_G(M), 1).
\eeqn
\item\beqn
\lker_{D(\kk G)}(\kk^{(G/M)_l})=\kk^{G}\bwt\mtr{core}_G(M)
\eeqn
\item\beqn
\lker_{\kk G}(\kk^{(G/M)_l})=\mtr{core}_G(M)
\eeqn\ene
\ep
\bpf
\onh note that $\kk^{(G/M)_l}\in \rep(G)$ and by Brauer's theorem for $\rep(G)$ gives that $\;_{D(\kk G)}<\kk^{(G/M)_l}>=\rep(G/\core_{G}(M))$.
\epf
Thus $\cd(\kk^{(G/M)_{l}})'=\cs(\mtr{core}_G(M), \{1\}, 1)'=\cs(\{1\}, \mtr{core}_G(M), 1)=\rep(G/\mtr{core}_G(M))=<\kk^{(G/M)_l}>$
 Thus previous conjecture is verified for $A=\kk^{G}$.

\section{On the Drinfeld double $D(H_8)$}
We show that for the unique eight dimensional semisimple noncommutative and noncocommutative Hopf algebra $H_8$ Conjecture \ref{coidconj} holds.
 \section{Drinfeld double of the eight dimensional Hopf algebra}
 
\subsection{Presentation of $H_{8}$}
Let $H_8$ be the unique eight dimensional noncommutative and noncoccomutative Hopf algebra \cite{}. Then $H_8$ is a semisimple Hopf algebra presented by generators $x, y, z$ and  the following relations:
\bq
x^2=y^2=z^2=1, xz=zx, zy=yz, xy=yxz
\eq

Comultiplication in $H_{8}$ is given by
\beq\label{x}
\D(x)=xe_{0}\ot x+ xe_{1}\ot y, \D(y)=ye_{0}\ot y+ye_{1}\ot x
\eeq
and $\D(z)=z\ot z$. the antipode is given by 
\beq\label{x}
S(x)=xe_{0}+ye_{1}, S(y)=ye_{0}+xe_{1}, 
\eeq
and $S(z)=z$. Here $e_{0}=\frac{1}{2}(1+z)$ and $e_{1}=\frac{1}{2}(1-z)$ are central idempotents of $A$.
\md
All normal Hopf subalgebras of $H_8$ are $k, K_1, K_2, H_8$ where $K_1=K(H_8)\cong k\Z_2$ and $K_2\cong k\Z_2\times k\Z_2$ coincides to the group algebra of group-like elements of $H_8$.
%\md All the characters here are central since $H_8$ is generated by two group like elements $x, y$ and a central element $z$
\subsection{Simple subcoalgebra of dimension $2$}
It follows that
$c_{11}=xe_{0}, c_{12}=xe_{1}, c_{21}=ye_{1},
 c_{22}=ye_{0}$.
 
 \subsection{Linear characters of $H_8$}
%With the notations from \cite{bd} let $g_i:=\hat{u}_i$. 
Let $u_1$, $u_2$ be the linear characters of $H_8$ that are not central. Then $u:=u_1u_2$ is a central linear character of $H_8$. Using the notations from \cite{bd} and \cite{aloui} one has that
$u_i(z)=1$ \md and $u_1(x)=-1$, $u_1(y)=1$  \md and $u_2(x)=1$, $u_2(y)=-1$.
\md
It follows that
\md $u(x)=u(y)=-1$ and $u(z)=1$.
\md Since $8p_{1}$ is the regular character of $H_{8}$ one has
\beq
\ch(x)=\frac{-(\eps+u_{1}+u_{2}+u)(x)}{2}=0
\eeq
Similalry $\ch(y)=\ch(xy)=0$ and $\ch(z)=-2$. Moreover $\ch(xz)=-\ch(x)=0=\ch(y)$.
Then the idempotent associated to $u_1u_2$ is $\eps_2$ while the idempotents of $u_1, u_2$ are $\eps_3$ and $\eps_4$ respectively.
\subsection{Central idempotents of $H_{8}$}
$e_{1}$ is the central idempotents of $\ch$.
\subsubsection{The 2 dimensional $H_{8}$-module}
Ohe has $M=kv_{1}\oplus kv_{2}$ with structure given by

$x.v_{1}=v_{2}$ and  $x.v_{2}=v_{1}$

$y.v_{1}=iv_{2}$ and $y.v_{2}=iv_{1}$

$zv_{1}=-v_{1}$ and $zv_{2}=-v_{2}$.
\subsection{Simple comodules associated to $H_{8}$}
One has $G(H_{8})=\Z_{2}\times\Z_{2}$ with $1, z=g_{1}g_{2}$ central grouplike elements and $g_{1}=xy(e_0 + ie_1)$ and $g_{2}=xy(e_0 - ie_1)$.
\subsection{Notations}

With the notations fom \cite{bd} it follows that $g_{i}=\widehat{u_{i}}$ and $z=\widehat{u_{1}}\widehat{u_{2}}$
\subsection{The left adjoint action} Following \cite{alaui} one has
\beq
x.x=x,\;x.y=yz,\;x.(yz)=y, \;x.(xz)=xz, x.(xy)=xyz=yx
\eeq	
and 
\beq
y.x=xz,\;y.y=y,\;y.(yz)=yz, \;y.(xz)=x, y.(yx)=yxz
\eeq	

It follows that \md $x.(xe_{0})=xe_{0}$, $x.(xe_{1})=xe_{1}$, $x.(ye_{0})=ye_{0}$, $x.(ye_{1})=-ye_{1}$ \md and
\md
$y.(xe_{0})=xe_{0}$, $y.(xe_{1})=-xe_{1}$, $y.(ye_{0})=ye_{0}$, $y.(ye_{1})=ye_{1}$
\subsection{Another basis for the center}

From here it follows that $1, z, xe_{0}, ye_{0}$ are central elements of $H_{8}$. Moreover also $xy+yx $ is central.

\subsection{$H_8$ regarded as a $D(H_8)$-module.} Using the notations from \cite{aloui} one has that
$
I(1)=H_8=W_1\oplus W_2\oplus W_3\oplus M_{1}\oplus M_{2}
$
where\md
$W_1=k1 \cong V_{1, \eps}$\md
$W_2=kz \cong V_{z, \eps}$\md
$ W_3=kg_{1}\oplus kg_{2}\cong V_{9}$\md
$M_{1}=kye_{0}\oplus kxe_{1}$\md
$M_{2}=kxe_{0}\oplus kye_{1}$.
%and $M_2=<xe_0, ye_1>$.
\subsection{Identification from table for $M_{{1}}$ and $M_{2}$}
 Since $M_{2}\dw_{H_{8}} \cong V_1\dw_{H_{8}}\cong \eps+u_{1}$ it follows that $M_{2}\cong V_{1}$ from \cite{bd}.   Similalry
  $M_1=<xe_1, ye_0> \cong V_2$ since $M_{1}\dw_{H_{8}} \cong V_2\dw_{H_{8}}\cong \eps+u_{2}$.

\subsection{Left kernels of modules} It can be computed that
\bne
\item The left kernel \beq
L_1:=\lker_{H_8}(u_1)=k1\oplus kz \oplus kye_0\oplus kxe_1=V_{1, \eps}\oplus V_{z, \eps}\oplus M_{1}.\eeq
where $M_{1}:=kye_{0}\oplus kxe_{1}$.
\item The left kernel 
\beq
L_2:=\lker_{H_8}(u_2)=k1\oplus kz \oplus kxe_0\oplus kye_1=k1\oplus kz\oplus M_{2}
\eeq
where $M_{2}:=kxe_{0}\oplus kye_{1}$.
\ene
\subsection{Fusion subcategory generated by $L_{1}$}
Note that
\beq
M_{1}^{2}=W_{1}\oplus W_{2}\oplus V_{g_{1}, u_{1}}\oplus V_{g_{2}, u_{1}}
\eeq
It follows that the fusion subcategory generated by $L_{1}$ is eight dimensional with objects $\{W_{1},W_{2},V_{g_{1}, u_{2}}V_{g_{2}, u_{2}}, M_{1}\}$.
\subsection{Fusion subcategory generated by $L_{2}$}
Note that
\beq
M_{2}^{2}=W_{1}\oplus W_{2}\oplus V_{g_{1}, u_{1}}\oplus V_{g_{2}, u_{1}}
\eeq
It follows that the fusion subcategory generated by $L_{2}$ is eight dimensional with objects $\{X_{1, \eps},W_{2},V_{g_{1}, u_{1}}V_{g_{2}, u_{1}}, M_{2}\}$.

\subsection{$L_{2}$ as an algebra.}
Note that $L_{2}$ can be presented as 
\beq
L_{2}=<1, z,b, c\:|\;bz=zb=b, cz=zc=-c, bc=cb=0, b^{2}=e_{0}, c^{2}=e_{1}>
\eeq
where $b=xe_{0}$ and $c=ye_{1}$
It is a commutative algebra with liner characters, $\eps, \nu=u_{1}|_{L_{2}}, \eta_{1}, \eta_{2}$ which are given as following:
$\eps: (1, z, a,b)\mapsto (1, 1,-1,0)$
\md
$\nu: (1, z, a,b)\mapsto (1, 1,1,0)$
\md
$\eta_{1}: (1, z, a,b)\mapsto (1, -1,0,1)$
\md
$\eta_{2}: (1, z, a,b)\mapsto (1, -1,0,-1)$.
Thus one has that 
$\eps|_{L_{2}}=\eps$
\md
$u_{1}|_{L_{2}}=\nu$
\md
$u_{2}|_{L_{2}}=\eps$
\md
${u_{1}|}_{L_{2}}=\nu$
\md
$\ch|_{L_{2}}=\eta_{1}+\eta_{2}$.
\subsection{Fusion subcategory $\cd(L_{2})$}
It follows that $\cd(L_{2})$ is formed by all simple modules whose restriction o $H_{8}$ contains only $\eps$ and $u_{2}$.
From the table of \cite{bd} it follows that
\beq
\cd(L_{2})=\{V_{1,\eps}, V_{g_{1}g_{2}, \eps}, V_{g_{1}, u_{2}}, V_{g_{2}, u_{2}}, M_{1}\}
\eeq
\blue{Then the conjecture can be written as
\beq
\{V_{1,\eps}, X_{z, \eps}, X_{g_{1}, u_{2}}, X_{g_{2}, u_{2}}, M_{1}\}'=\{X_{1,\eps}, X_{z, \eps},X_{g_{1}, u_{1}}, X_{g_{2}, u_{1}}, M_{2}\}.
\eeq}
\bab which is equivalent to $M_{1}$ centralizes $M_{2}$.\eab
\subsection{Remark on normal coidela subalgebras}
All normal coideal subalgebras of $H_{8}$ are commutative.
\subsection{Normal Hopf subalgebrs of $D(H_8)$}
Let $$k \subset K_1   \subset K_2 \subset H_8$$ be the lattice of normal Hopf subalgebras of $H_8$. One has $K_1=K(H_8)=\mathbb{C}1 \oplus \mathbb{C}z$ and $\dim_{\kk} K_1=2$ while $\dim_{\kk} K_2=4$ since $K_{2}=k1\oplus kz\oplus kg_{1}\oplus kg_{2}$. Moreover $K_{1}=k\Z_{2}$ and $K_{2}=k[\Z_{2}\times \Z_{2}]$.
\md
By duality, since $H_{8}$ is self dual it follows that $$k \subset P_{1}   \subset P_{2} \subset H_8^*$$ is the lattice of normal Hopf subalgebras of $H_8^*$ with $P_1:=(H_8//K_2) ^{*}$ and $P_2:=(H_8//K_1)^* $. Clearly $K(H_8^*)=P_1$.
\subsection{$D(H_{8})$ as an abelian extension} Follwing \cite{nrm} one has that
\beq
k \ra B(K_{1}, K_{1})\ra D(H_{8})\ra D(K_{1}, K_{1})\ra k
\eeq
\beq
k \ra B(K_{2}, K_{2})\ra D(H_{8})\ra D(K_{2}, K_{2})\ra k
\eeq
are exact abelian extensions. All the Hopf subalgebras have dimension $8$. THhus the quotients have also dimension $8$.
\subsection{Fusion subcategory generated by $V_{9}$}
Note that $V_{9}^{*}=V_{9}$ just by inspecting the table. Note also that $V_{9}^{2}$ has restriction $(2\eps+2u, 2+2z)$ to $H_{8}$ and $H_{8}^{*}$.
By inspecting the table this implies that $V_{9}^{2}$ cannot contain $2$-[dimensional representations. Thus
\beq
V_{9}^{2}=X_{1, \eps}\oplus X_{z, \eps}\oplus X_{u, 1}\oplus X_{u,z}
\eeq

\subsection{Fusion subcategories $\cd(K_{i})$ and $\cd(L_{i})$}
Since $K_{1}=K(H_{8})$ it follows that $\cd(K_{1})$ coincides to $F^{-1}(\rep(H_{8})_{ad})$. From the Table 1.1 of \cite{bd} one has that
\beqn
\co(\cd(K_1))=\{X_1, X_2,X_{3}, X_{4}, X_{5}, X_{6}, X_{7}, X_8, V_1, V_2, V_3, V_4, V_9, V_{10}\}=N_{1}
\eeqn
Also note from the same table
\beqn
\co(\cd(K_2))=\{ X_{1, \eps},  X_{1, u}, X_{z, u}, X_{z, \eps}, V_9\}
\eeqn
\onh $<K_{1}>=\{V_{1, \eps}, V_{z, \eps}\}$ thus one obtains that

By duality one has that 
\beqn
\co(\cd(P_{1}))=\{X_1, X_2,X_{3}, X_{4}, X_{5}, X_{6}, X_{7}, X_8, V_5, V_6, V_7, V_8, V_9, V_{10}\}
\eeqn
and 
\beqn
\co(\cd(P_2))=\{ X_{1, \eps},  X_{1, u}, X_{z, u}, X_{z, \eps}, V_{10}\}
\eeqn
\subsection{Normal Hopf subalgebras as $D(H_{8})$ modules}
One has $K_{1}=X_{1, \eps}\oplus X_{z, \eps}$ and  \md
$K_{2}=X_{1, \eps}\oplus X_{z, \eps}\oplus W_{3}$.\md Recall that $W_{3}=kg_{1}\oplus kg_{2} \cong V_{9}$.
\subsection{Fusion subcategories generated by normal Hopf algebras}
One has that 
$<K_{1}>=\{X_{1, \eps}, X_{z, \eps}\}$  and $<K_{2}>=\{X_{1,\eps}, X_{1, u}, X_{z, \eps}, X_{z,u}, V_{9}\}$.
By duality one obtains that
$<P_{1}>=\{X_{1, \eps}, X_{1, u}\}$ and 
$<P_{2}>=\{X_{1,\eps}, X_{1, u}, X_{z, \eps}, X_{z,u}, V_{10}\}
$.
\subsection{Centralizers of fusion subcategories $\cd(K_{i})$}
\onh one has that equality $\cd(K_{1})'=<K_{1}>$ gives that \beq\label{r1}
\{X_1, X_2, X_{3}, X_{4}, X_{5}, X_{6}, X_{7}, X_8, V_1, V_2, V_3, V_4, V_9, V_{10}\}'=\{X_{1, \eps}, X_{z, \eps}\}
\eeq
One has that equality $\cd(K_{2})'=<K_{2}>$ gives that 
\beq\label{r2}
\{ X_{1, \eps},  X_{1, u}, X_{z, u}, X_{z, \eps}, V_9\}'=\{X_{1,\eps}, X_{1, u}, X_{z, \eps}, X_{z,u}, V_{9}\}
\eeq
i.e $\cd(K_{2})'=\cd(K_{2})$.
\md
By duality one has that equality $\cd(P_{1})'=<P_{1}>$ gives that \beq\label{r3}
\{X_1, X_2, X_{3}, X_{4}, X_{5}, X_{6}, X_{7}, X_8, V_5, V_6, V_7, V_8, V_9, V_{10}\}'=\{X_{1, \eps}, X_{1, u}\}
\eeq
and equality $\cd(P_{2})'=<P_{2}>$ gives that
\beq\label{r4}
\{X_{1, \eps},  X_{1, u}, X_{z, u}, X_{z, \eps}, V_{10}\}'=\{X_{1,\eps}, X_{1, u}, X_{z, \eps}, X_{z,u}, V_{10}\}
\eeq
i.e. $\cd(P_{2})'=\cd(P_{2})$ 
\subsection{On the adjoint action}
By \cite{bd} it follows that $\rep(D(H_{8}))_{ad}$ has as simple objects $\{X_{1}, \cdots, X_{8}\}$. Since $\rep(D(H_{8}))$ is a modular tensor category one has that $\rep(D(H_{8})_{ad})'=\rep(D(H_{8}))_{pt}=\{X_{1}, \cdots, X_{8}\}$.
Thus
\beq\label{r5}
\{X_{1}, \cdots, X_{8}\}'=\{X_{1}, \cdots, X_{8}\}
\eeq
\subsection{Centralizer of normal fusion categories}
Note that 
\beqn
\cd(K_{1}, K_{1})=\cd(P_{2})\cap \cd(K_{1})=\cd(P_{2})=\{ X_{1, \eps},  X_{1, u}, X_{z, u}, X_{z, \eps}, V_{10}\}
\eeqn
\onh
\beqn
\cd(K_{2}, K_{2})=\cd(P_{1})\cap \cd(K_{2})=\cd(K_{2}).
\eeqn
\md Note that
$B(K_{1}, K_{2})$ is a normal Hopf subalgebra. One has that 
$\cd(K_{1}, K_{2})=\cd(P_{2})\cap \cd(K_{2})=\{ X_{1, \eps},  X_{1, u}, X_{z, u}, X_{z, \eps}\}$. On the other hand
$\cd(K_{2}, K_{1})=\cd(P_{1})\cap \cd(K_{1})=\{ X_{1},\cdots  X_{8}, V_{9}, V_{10}\}$. Thus
\beq\label{r6}
\{ X_{1, \eps},  X_{1, u}, X_{z, u}, X_{z, \eps}\}'=\{ X_{1},\cdots,X_{8}, V_{9}, V_{10}\}.
\eeq

Similalrly $B(K_{2}, K_{1})$ is a normal Hopf subalgebra and 
$\cd(K_{1}, K_{2})=\cd(P_{2})\cap \cd(K_{1})=\cd(P_{2})$. On the other hand $\cd(K_{2}, K_{1})=\cd(P_{2})\cap \cd(K_{1})=\cd(P_{2})$. Thus we get the same relation.

Also $B(H_{8}, K_{1})=K_{1}$ is normal and $\cd(H_{8}, K_{1})=\cd(K_{1})$. On the other hand $\cd(K_{1}, H_{8})=\cd(H_{8})$.
\subsection{Left kernels of normal Hopf algebras}
One has $\lker_{H_{8}}(K_{1})=H_{8}$ and by Brauer's theorem $\lker_{D(H_{8})}(K_{1})=(A//K_{1})^{*}\bwt H_{8}$. Similalrly
$\lker_{H_{8}}(K_{2})=K_{2}$ and by Brauer's theorem $\lker_{D(H_{8}}(K_{2})=(A//K_{2})^{*}\bwt K_{2}$. The conjecture on left kernels is satisfied.
\subsection{Left kernels of coideal subalgebras}

By Brauer's theorem one has that $\dim_{\kk} \lker_{D(H_{8})}(L_{2})=8$. On the other hand $\lker_{H_{8}}(L_{2}|_{H_{8}})=\lker_{H_{8}}(u_{1})=L_{1}$. Moreover $\lker_{H_{8}^{*}}(L_{2}|_{{H_{8}^{*}}})=\lker_{H_{8}^{*}}(d)=k$.\mdn
\blue{It verifies that $\lker_{A}(\lker_{A}(L_{i}))=\lker_{A}(L_{i})$. It is also true for the normal Hopf subalgebras.}
\subsection{Stable characters}
One has that $K_{1}=k1\oplus kz$ has both characters stable.
On the other hand $K_{2}=k1\oplus kz\oplus kg_{1}\oplus kg_{2}$ has to stable characters, the trivial character and $\eta: (1, z, g_{2}, g_{2})\mapsto (1, 1, -1,-1)$. Moreover $u_{i}\dw_{K_{2}}=\eta$.
\subsection{Finish description of normal Hopf subalgebras of $D(H_{8})$}
One has $B(K, L, G)$ with $G\subset G_{A}(K)\cap G_{A^{*}}((A//L)^{*})$. In our situation
 $G\cong \Z_{2}$ and the dimension of the Hopf algebra doubles. It should be of dimension $16$ and therefore need to determine when the one dimensional characters are central in $D(H_{8})^{*}$.
 \md
 Since the character of $V_{\eta, g}$ is $\eta\ot g$ it follows that $V_{1, \eps}, V_{z, \eps}, V_{z, u}, V_{1, u}$ are central and except the first one (which is the trivial module) all the other have kernels of dimensions $32$.
 \md
 There are also two normal Hopf subalgebras of dimension $16$ as the kernels of $V_{z, u_{1}}+V_{z, u_{2}}$ and $V_{g_{1}, u}+V_{g_{2}, u}$.
 \subsection{Complete the other Hopf algebras}
 \bne
 \item
 $
 B(K_{1}, K_{1}, G)=P_{2}\bwt K_{1}\oplus (C_{u_{1}}\oplus C_{u_{2}})\bwt C_{d}
 $
 \md
$ \dim_{\kk}  B(K_{1}, K_{1}, G)=16$. Since $\cd(K_{1}, K_{1})=\cd(P_{2})\cap \cd(K_{1})=\cd(P_{2})=\{ X_{1, \eps},  X_{1, u}, X_{z, u}, X_{z, \eps}, V_{10}\}$
it follows that \md
$
\cd(K_{1}, K_{2}, G)=\{ X_{1, \eps},  X_{1, u}, X_{z, u}, X_{z, \eps}\}.
$
Thus $V_{10}$ is central.
\item
 $
 B(K_{1}, K_{2}, G)=P_{2}\bwt K_{2}\oplus C_{\ch}\bwt C_{d}
 $
 it has dimension $32$.
 \item
 $
 B(K_{2}, K_{1}, G)=P_{1}\bwt K_{1}\oplus (C_{u_{1}}\oplus C_{u_{2}})\bwt (kg_{1} \oplus kg_{2})
 $
 \md
 \ene
 \blue{Is $D(H_{8})$ an indecomposable fusion category? Is the centralizer of any normal also normal.}
 \subsection{Stable characters of $K_{1}$}
 One has two classes $(H_{8})_{ad}$ seating over $\eps$ and $\ch$ seating over the other character.
 \subsection{Stable characters of $K_{2}$}
 $\hker_{H_{8}}(u)=K_{2}$
 One has that $u_{1}$ and $u_{2}$ are seating over the same stable character.
 \subsection{Cosets generated by the adjoint fusion subcategory}
 \bne
  \item
 $\cd_{0}=\{X_{1}, \cdots, X_{8}\}$
 \item
 $\cd_{1}=\{V_{9},\; V_{{10}}\}$
  \item
 $\cd_{2}=\{V_5,\; V_7\}$
  \item
 $\cd_{3}=\{V_8 ,\;V_6\}$
 \item
 $\cd_{4}=\{V_1,\; V_4\}$
  \item
 $\cd_{5}=\{\}$
  \item
$ \cd_{6}=\{\}$
  \item
$ \cd_{7}=\{\}$
  \item
 $\cd_{8}=\{\}$
 \ene
 \subsection{Normal union subcategories generated by equivalence classes of adjoint cosets}
 One has $8$ normal fusion subcategories of dimension 16 $N_{i}=\cd_{0}\vee \cd_{i}$. They correspond to Hopf subalgebras of dimension 4.
 \bne
 \item
 $ N_{1}=\cd_{0} \vee \cd_{1}=\{\cd_{0}, V_{9}, V_{10}\}=\cd(K_{2}, K_{1})=\rep(D//B(K_{2}\bwt K_{1}))$
 \md Is this Hopf subalgebra central? yes since it is a group algebra
 \md
 $N_{1}'=\{X_{1, \eps}, X_{1, u}, X_{z, \eps}, X_{x, u}\}$
  \item 
 $ N_{2}=\cd_{0} \vee \cd_{2}=\{\cd_{0}, V_{5}, V_{7}\}$
  \item
   $ N_{3}=\cd_{0} \vee \cd_{3}$ 
  \item 
   $ N_{4}=\cd_{0} \vee \cd_{4}$
  \item 
   $ N_{5}=\cd_{0} \vee \cd_{5}$
   \item 
    $ N_{6}=\cd_{0} \vee \cd_{6}$
    \item
    $ N_{7}=\cd_{0} \vee \cd_{7}$
 \item
 $ N_{7}=\cd_{0} \vee \cd_{7}$
  \ene
 \md
 
 It can be checked that any $16$ dimensional  normal fusion subcategory is of this type.
 \md
 It is the centralizer of a four dimensional fusion subcategory.
 \subsubsection{Central Hopf subalgebras of dimension 4}
 Consider generators for $\ovr{G}(D(H_{8}))$ 
 
 $a=\eps\bwt z$, 
 
 $b=g_{1}\bwt u_{1}$ 
 
 and 
 
 $c=g_{2}\bwt u_{2}$.

 Compute
 $\rep(D/\ovr{G}(a, b))$
They correspond to this central Hopf subalgebras of dimension $4$.
 \subsection{Fusion subcategories of dimension 32}
 \subsubsection{Central  Hopf subalgebras of dimension $2$}
 \beq
 \rep (D(A)/\ovr{G}(a))=\{\cd_{0}, \cd_{1}, \cd_{2}, \cd_{3},\cd_{1}\cd_{2}\}
 \eeq
 \begin{tabular}{ l | c r r r r}
   & $\eps$ & $u_{1}$ & $u_{2}$ & u & $\ch$ \\  \hline
            1     & 1 & 1 & 1 & 1 & 2\\
  $x$ & 1 &   -1 & 1 & 1 & 0\\
  $y$ & 1 &    1 & -1 & 1 & 0\\
  $z$ & 1 &    1 & 1 & 1 &-2
\end{tabular}

and

\begin{tabular}{ l | c r r r}
   & $\eps$ & $u_{1}$ & $u_{2}$ & u\\
  \hline
  1 & 1 & 1 & 1 & 1\\
  $g_{1}$ & 1 & -1 & -1 & 1\\
  $g_{2}$ & 1 & -1 & -1 & 1
  \\ g & 1 & 1 & 1 & 1
\end{tabular}
\newpage
\bibliographystyle{amsplain}
\bibliography{mec-rep-drd}
\section{Tambara-Yamagami category}
\subsection{Description of simple objects}
$Z_{\ro}$ where $\ro$ has to satisfy the identity
\beqn
\ro(ab)=\ch(a,b)^{-1}\ro(a)\ro(b)
\eeqn
There are four possibilities of such $\ro's$ sending $(b,c)\mapsto (\pm i, \pm i)$. Using the fact $\Delta^{2}=1. i$ one obtains simple objects whose restrictions to $\cc$ is $m$.

Simple $Z$ objects are $Z_{\pm i, \pm i\;, \pm \sqrt{\Delta}}$ where $\Delta^{2}=1$ if the values of $\ro$ in $a, b$ coincide and $\Delta^{2}=\pm i$ otherwise.
\br
$V_{7}, V_{8}$ are $Y's $ for the dual category oof representations.
\er
\subsubsection{Settings}

Consider $A=\Z_2 \times \Z_2=\{ 1,a, b, c\}$. Consider $c= u_ 1 u _ 2$, $a = u_ 1$, $b= u _ 2$.

One has that
\beq
\ch_{c}(a,a)=\ch_{c}(b,b)=-1, \; \ch_{c}(a, b)=1
\eeq

It follows that 
\beqn
\ch ( c, a) = \ch (b, a)\ch(a, a) = -1=
\ch(c, b)\eeqn
and 
\beqn
\ch ( c, c) Ê= \ch (b, c) \ch (a, c) = 1.
\eeqn

\subsection{Centralizers of invertible objects}
 $X_{a,\eps}$ centralizes any $X_{b, c} $.

It follows that $X_{a,\eps}$ centralizes 
$Y_{b, c}$ if and only if $\ch(a, bc)= 1$

\bne
\item $a=1$ all of them in the centralizers
\item $b= -c$ are in the centralizers
\ene

It follows that $X_{a,\eps}$ centralizes 
$Z_{\ro, \delta}$ with $\eps = 1$ and $\ro( a) =1$.

\subsection{Centralizers of $Y_{ b, c}$}
\bne
\item
$Y_{b, c}$ does not centralize any $Z_{\ro, \delta}.$
\item
$X_{a, \eps}$ centralize $Z_{ \ro, \delta}$ if and only if $\eps\ro(a)=1$.
\item
$Y_{l,m} $Êcentralizes $Y_{ b, c}$ if and only if
$\ch(m, b)=\ch(l, c)$ and $\ch(m, c)= \ch (l,b)$
\item
 $Z_{\ro, \delta}$ centralizes $Z_{ \ro, \delta'}$
if and only if $\ro(a)\ro(a')=\Delta\Delta'$ for all $a, a' \in A$.
\ene
\subsubsection{
Recognizing the indexed modules}
\bne
\item $X _{ 1, 1}= ÊX_{1, \eps}$
\item $X _ {1 -1}= X_{z, \eps}$
\item $X_{a, i} = X_{u_1, g_1}$
\item $X_{a, -i}=X_{ u_1, g_2}$
\item
 $X_{b, i}=X_{u_2, g_{1}}$
 \item
 $X_{b, -i}=X_{u_{2}, g_{2}}$
  \item
 $X_{c, 1}=X_{u, 1}$
 \item
 $X_{c, -1}=X_{u,z}$
\item $V$-urile sunt $Y_{a+b}$, where $a+b$ is the restriction to $A$.
\subsection{Writing down the centralizer}
\beqn
(X_{1, -1})'=\{X_{1}, \cdots X_{8}, Y_{a, c}, Y_{1, b}, Y_{1, c} , Y_{b,c}, Y_{a,b}, Y_{1,a}\}
\eeqn
\beqn
(X_{a, i})'=\{X_{1}, \cdots X_{8}, Y_{a, c}, Y_{1, b}, Z_{\ro(a)=-i}\}
\eeqn
\beqn
(X_{a, -i})'=\{X_{1}, \cdots X_{8}, Y_{a, c}, Y_{1, b}, Z_{\ro(a)=+i}\}
\eeqn
\beqn
(X_{b, -i})'=\{X_{1}, \cdots X_{8}, Y_{b, c}, Y_{1, a}, Z_{\ro(b)=+i}\}
\eeqn
\beqn
(X_{b, i})'=\{X_{1}, \cdots X_{8}, Y_{b, c}, Y_{1, a}, Z_{\ro(b)=-i}\}
\eeqn
\beqn
(X_{c, 1})'=\{X_{1}, \cdots X_{8}, Y_{b, a}, Y_{1,c}, Z_{\ro(c)=1}\}
\eeqn
\beqn
(X_{c, -1})'=\{X_{1}, \cdots X_{8}, Y_{1, c}, Y_{a,b}, Z_{\ro(c)=-1}\}
\eeqn
\beqn
(Y_{1, a})'=\{X_{1, 1},X_{1, -1}, X_{b, -i}, X_{b, +i}, Y_{1,b} \}
\eeqn
\beqn
(Y_{1, b})'=\{X_{1, 1}, X_{1, -1}, X_{a, -i}, X_{a, +i}, Y_{1,a} \}
\eeqn
\beqn
(Y_{1, c})'=\{X_{1, 1}, X_{1, -1}, X_{c, 1}, X_{c, -1}, Y_{1,c} \}
\eeqn
\beqn
(Y_{a, b})'=\{X_{1, 1}, X_{1, -1}, X_{c, 1}, X_{c, -1}, Y_{a,b} \}
\eeqn
\beqn
(Y_{b, c})'=\{X_{1, 1}, X_{1, -1}, X_{b, i}, X_{b, -i}, Y_{a,c} \}
\eeqn
\beqn
(Y_{a, c})'=\{X_{1, 1}, X_{a, i}, X_{b, i}, Y_{b,c} \}
\eeqn
\beqn
(Z_{-i, -i, \omega})'=\{X_{1, 1}, X_{1, -1}, X_{a, i}, X_{a, -i}, Z_{} \}
\eeqn
\ene

\newpage
\blue{It uses the dual category and therefore cannot be generalized to any fusion category in a straightforward manner. It would be interesting to find a correspondent approach for any fusion category.}

\bb{
Schneider's map for Drinfeld double $G(A^*)\ra G(A)\cap Z(A)$ is given by both restrictions
\beqn
V_{g, \eta} 
\eeqn
iff $\eta\bwt g$ is central. %\onh is a twits.
Maybe this is an equivalence relation for factorizable Drinfeld doubles
}

\newpage
\bn{cor} \label{onechr'}Let $A$ be a factorizable Hopf algebra. Let $\ch \in \Irr(A)$ be an irreducible character of $A$ and $E_{\ch}:=\phi_{R}^{-1}(e_{\ch}) \in C(A)$ be the associated primitive idempotent of $\ch$. Then $$<\ch>'=\Rep(A//\prod_{\{d\;|\;E_{\ch}(d)\neq 0\}}L(d)).$$
\end{cor}
\blue{It is kind of a duality?? all $\ch_{j}$ such that $E_{j}(d)\neq 0$ for all $d $ with $E_{\ch}(d) \neq 0$. what can use the bijection  in this form}

\bpf
By Theorem \ref{interms} and Equation \eqref{centrofint} one has that 
\
\begin{eqnarray*}
 \Rep(A//\prod_{\{d\;|\;E_{\ch}(d)\neq 0\}}L(d))'&=&(\bigvee_{d\;|\;E_{\ch}(d)\neq 0}\Rep(A//L(d)))'\\ &=&\bigcap_{d\;|\;E_{\ch}(d)\neq 0}\Rep(A//L(d))'.
\end{eqnarray*}
\epf
%\bn{rem}\lb{centrofz}Something similar holds for $Z_{ _d}$. Namely $\ch_j\in Z_d'$ if and only if $E_j(d)\neq 0$ and $s_{ij}=\omega_id_id_j$. The scalar $\omega_i$ does not depend on $j$. %It follows that $\ch_i \in ker_{A^*}(d)$ if and only if $\ch_i$ centralize each $\ch_j$ with the property that $E_j(d)\neq 0$.
%\er
%\blue{If $A$ is just quasitriangular does a similar formula holds? Need to rewrite the centralizer; but a formula does not exists since $C(A)$ not bijectively to $Z(A)$.
%}
\newpage
\bp\label{inters}
Let $L$ and $K$ be two normal Hopf subalgebras of a semsimple Hopf algebra $A$. Then 
\beqn LK//L\cap K \cong (L//L\cap K) \otimes (K//L\cap K)\eeqn 
as Hopf algebras.
\ep
\bpf Using previous Theorem it follows that 
\beqn 
\dim_{\kk}(L\cap K)=\frac{\dim_{\kk} L \;\dim_{\kk} K}{\dim_{\kk} LK}.
\eeqn Applying the same formula for $L//L \cap K$ and  $K//L \cap K$ inside $A//L \cap K$ it follows that  the two normal Hopf subalgebras of $ A//L\cap K$ have trivial intersection Then by \cite[Proposition 2.4]{cejm} it follows that $LK//L\cap K$ is the tensor product of the two Hopf algebras.
\epf

\newpage
\label{repalg}
Let $A$ be a semisimple Hopf algebra and $V_1=k$, $V_2,\cdots, V_s$ be the homogeneous $D(A)$-submodules of $A$, thus
\beqn
A=V_{1}\oplus V_{2}\oplus\cdots \oplus,  V_{s}.
\eeqn

The following Lemma follows from the proof \cite[Theorem 6.3]{KSZ}. For the sake of completeness we sketch its proof below.

\bl\cite{KSZ}\lb{sommyrep}
Let $V_i$ as above be a homogenous component of the $D(A)$-module $A$. Then $H^*S(p_i) \cong V_i$ as $D(A)$-modules.
\el

\bpf
It was proven in \cite[Theorem 5.13]{repalg} that $p_V$ is a central idempotent of $C(A)$. Following \cite{Zind} one has that $\End_{\Dr(A)}(A)=C(A)$ and the isomorphism is given by $\ch \mapsto (\; \lh \ch)$.
\epf
\newpage
\bpf
Note that any primitive idempotent of $A^*$ is of the type $\blue{X^{d}_{pp}}$ for a coalgebra matrix $\{x^{d}_{pq}\}$.
Writing $E_j$ as sum of primitive idempotents of $A^*$ one obtains
that
\beq
E_{j}=\sum_{d,p \in \ca_{j}}X^{d}_{pp}
\eeq
 and the above inequality follows since $X^{d}_{pp}(d')=\delta_{d,d'}$. Since $\sum_{j=1}^{s}E_{j}=\eps$ inside $C(A)$ the above equality also follows.\epf
In this situation the second relation of Corollary \ref{norm'} becomes:

\begin{eqnarray*}
 \Rep(A//K)'=<\ch_V\;|\;E_{V\ot V^*}(x)\neq 0 \;\text{for some}\; x \in \Irr(K^*)>.
\end{eqnarray*}
\bb{This needs the formula restriction for factorizable}
Need a description of $E_{V\ot V^{*}}$.
\blu{ Study the induced functor by $\delta_R$.
It can be verified that if the conditions of the above theorem are satisfied then the isomorphism $\delta_R$ is an isomorphism of quasitriangular Hopf algebras.}
\subsection{The maps $\pi_R$ and $\pi_{\tilde{R}}$} Let $A$ be a quasitriangular Hopf algebra and consider the map $\pi_R:D(A) \ra A$ given by $f\bowtie a \mapsto ((f\ot id)(R))a$ surjective Hopf algebra map. Moreover $\pi_R$ is a retraction of the inclusion map.
 
This map appears in [\cite{dr}, Proof of Proposition 6.2] where it was observed that $\pi_R$ is a Hopf algebra homomorphism. Since  for $\tilde{R} := (R_{21})^{-1} = S^{-1}(R_2) \ot R_1 = R_2 \ot S(R_1)$ one has that $(A,\tilde{ R})$ is again a quasitriangular Hopf algebra \cite[VIII.2]{Kas}, there is also a projection $\pi_{\tilde{R}}: D(A) \ra A$. \blue{are these maps needed}

\subsection{On the quantum double of a factorizable Hopf algebra} Let $(A;R)$ be a finite-dimensional quasitriangular Hopf algebra.
Define
$$
F_R := 1\ot R_2 \ot R_1 \ot 1.
$$
Then $F_R$ is a 2-cocycle for the tensor product $A\ot A$
with componentwise bialgebra structure,

By \cite[Theorem 2.3]{schfact} note that for any quasitriangular Hopf algebra $A$ the map
\beqn
\delta_{R}: D(A) \ra (A\ot A)^{F_R}, \;\; x\mapsto (\pi_{\tilde{R}} \ot \pi_{R}\pi_{\tilde{R}})(\Delta_{D(A)}(x))
\eeqn 
is a Hopf algebra homomorphism. Moreover by loc.cit. $(A,R)$ is factorizable if and only if $\delta_{R}$ is an isomorphism.

\bp\lb{restrfact}
Under  the isomorphism $\delta$ from above the restriction of a left $D(A)$-module $M\boxtimes N$ to $A$ is given by $$(M\bx N)\dw^{D(A)}_A
\cong M\ot N.$$ 
\ep
\bpf
It is easy to check that $\delta_{R}(\eps \bwt a)=a_1\ot a_2$.
\epf

\blue{Here the maps $\pi$ are used but this proposition not needed if the other results are not included.}
% In this case one can take  $v=u_D=\sum_{i=1}^nS^{2}b_i^* \bowtie b_i$.

 \br
It follows from above that $\hker_A(M)$ is the largest Hopf subalgebra of $A$ contained in $\lker_A(M)$. In fact it is easy to see that $\hker_A(M)$ is the largest subcoalgebra of $A$ contained in $\lker_A(M)$. Indeed if $C\subset \lker_{A}(M)$ is a subcoalgebra of $A$ then $\sum_{n}C^{n}\subset \lker_A(M)$. On the other hand note that $\sum_{n}C^{n}$ is a Hopf subalgebra of $A$ and therefore one has that $\sum_{n}C^{n}\subseteq \hker_A(M)$. \blue{not really needed}
 \er
 \subsection{Centralizers from normal closures}
\subsection{Here}

For $K$ a Hopf subalgebra of $A$ denote by $\cD(K)$ to be the fusion subcategory of $\mtr{Rep}(D(A))$ generated by all modules having trivial $K$-action. Clearly $\cD(K)=\Rep(D(A)//L(K))$ where $L(K)$ is the smallest normal left coideal subalgebra of $D(A)$ containing $K$.
\section{Centralizers for fusion subcategories of Deligne products and their centralizers}
\subsection{Fusion subcategories of direct products of categories} First we need a preliminary result for the structure of fusion subcategories of a Deligne product of two fusion categories.\md 
Let $\mtc{C}^1$ and $\mtc{C}^2$ be two fusion categories. Identify
them with the corresponding fusion subcategories of $\mtc{C}^1\boxtimes\mtc{C}^2$ given by $\mtc{C}^1\boxtimes 1$ and $1 \boxtimes \mtc{C}^2$ respectively. Then every simple object of
$\mtc{C}^1\boxtimes\mtc{C}^2$ is of the form $X_1\boxtimes X_2$
where $X_i$ is a simple object of $\mtc{C}^i$.
%%%
Let $\mtc{D} \subset \mtc{C}^1\boxtimes\mtc{C}^2$ be a fusion
subcategory. Define $\mtc{L}^i(\mtc{D}):=\mtc{D}\cap \mtc{C}^i$,
$i=1,2$. Let also $\mtc{K}^1(\mtc{D})$ be the fusion subcategory
generated by all simple objects $X_1$ of $\mtc{C}^1$ such that $X_1
\boxtimes X_2 \in \mtc{D}$ for some simple object $X_2$ of
$\mtc{C}^2$. Similarly define the fusion subcategory
$\mtc{K}^2(\mtc{D})$. Clearly
$\mtc{L}^1(\mtc{D})\boxtimes\mtc{L}^2(\mtc{D})\subset \mtc{D}
\subset\mtc{K}^1(\mtc{D})\boxtimes\mtc{K}^1(\mtc{D})$.
The following Lemma is taken from \cite{dgno}.
\bt\label{dg}
Let $\mtc{D} \subset \mtc{C}^1\boxtimes\mtc{C}^2$ be a fusion subcategory. Then there is
a group $\mtc{X}$ and faithful $\mtc{X}$-gradings $\mtc{K}_i(\mtc{D})=\oplus_{x
\in \mtc{X}}\mtc{K}^i(\mtc{D})_x$ with trivial components
$\mtc{L}^i(\mtc{D})$, $i=1, 2$ such that
\begin{equation}\label{81}
    \mtc{D}=\oplus_{x \in \mtc{X}}\mtc{K}^1(\mtc{D})_x\boxtimes\mtc{K}^2(\mtc{D})_x.
\end{equation}
\et
\bt\lb{interchange}
Let $\cD$ be a fusion subcategory of a Deligne product of two modular categories $\cC_1$ and $\cC_2$. Then 
$\cK_i(\cD')=\cL_i(\cD)'$ and $\cL_i(\cD')=\cK_i(\cD)'.$
Moreover the group $\cX$ is the same for $\cd$ and $\cd'$ and there is the following relation between gradings:
\et

\bp
If $\cD$ is Lagrangian fusion subcategory of $\cc_{1}\bx  \cc_{2}$ then $K_i(\cD)'=L_i(\cD)$ and $L_i(\cD)'=K_i(\cD).$\ep

\bp Let $\cc_{1}$ and $\cc_{2}$ be two braided fusion categories and form $\cc_{1}\bx \cc_{2}$.
Let $\cD$ and $\cE$ be two fusion subcategories of $\cC_1\bx \cC_2$. Then one has the following equalities

1) $\cK_i(\cD\cap \cE)=\cK_i(\cD)\cap \cK_i(\E)$

2) $\cK_i(\cD\vee \cE)=\cK_i(\cD)\vee \cK_i(\E)$

3) $\cL_i(\cD\cap \cE)=\cL_i(\cD)\cap \cL_i(\E)$

4) $\cL_i(\cD\vee \cE)=\cL_i(\cD)\vee \cL_i(\E)$
\ep

\bpf
The proof of the first two items is straightforward. For the last two items use Theorem \ref{interchange}.
Indeed one has that
\begin{eqnarray*}
 \cL_i(\cD\cap \cE)'&=&\cK_i((\cD\cap \cE)')=\cK_i(\cD' \vee \cE')=\\&=&\cK_i(\cD' ) \vee \cK_i(\cE')=\\&=&\cL_i(\cD)'\vee \cL_i(\cE)'.
\end{eqnarray*}
 Thus $\cL_i(\cD\cap \cE)=\cL_i(\cD)\cap \cL_i(\E).$
Similarly it can be proven that $\cL_i(\cD\vee \cE)=\cL_i(\cD)\vee \cL_i(\E).$
\epf

It is easy to see the following fact which also follows from a more general result for modular categories, namely that $\cZ(\cC)\cong \cC\boxtimes \cC^{\mtr{rev}}$ as braided fusion categories.

% now I have a formula for the M\"uger's centralizer

\subsection{The centralizer of a normal fusion subcategory of a modular category}
\bt\lb{form for normal}
Let $A$ be a factorizable semisimple Hopf algebra and $K$ be a normal Hopf subalgebra of $A$. Then
\beq\lb{Muegers}
\Rep(A//K)'\subseteq \cK_1(\Rep(D(A)//(A//K)^*\bwt \LKER_A(K)))
\eeq
\et

\bpf
As above let $\cD(K)$ be the fusion subcategory of $\Rep(D(A))$ with the property that $K$ acts trivially on each of its objects. Using Proposition \ref{justonecomp} it follows that 
\beqn
\cd(K)'\subseteq \rep(D(A)//(A//K)^*\bwt \LKER_A(K))
\eeqn
Note that by Proposition \ref{restrfact} one has $\cL_1(\cd(K))=\rep(A//K)$.
Thus Theorem \ref{interchange} gives that
\beqn\rep(A//K)' \subseteq \cK_1((\Rep(D(A)//(A//K)^*\bwt \LKER_A(K))))\eeqn 
and the proof is complete.\epf

Let $K$ be a normal Hopf subalgebra of $A$ and $K_{st}:=\cap_{x \in \mtr{St}_A(K)}\lker_K(x)$.

\bp
Let $A$ be a factorizable semisimple Hopf algebra and $K$ be a normal Hopf subalgebra of $A$. Then
\beq 
\Rep(A//L(K_{st}))'\subseteq \cK_1(\Rep(D(A)//(A//K)^*\bwt \LKER_A(K)))
\eeq
\ep
\bpf
Using Theorem \ref{interchange} \beqn\cK_1(\cD(K))'=\cL_1(\Rep(D(A)//(A//K)^*\bwt \LKER_A(K))). \eeqn But using Proposition \ref{restrfact} it is easy to see that $\cK_1(\cD(K))=\Rep(A//L(K_{st}))$.
\epf
\blue{\bt\lb{genchr}
Let $\cD:=\cD(K,\;L,\;\cX,\;\psi)$ be a a normal fusion subcategory of $\Rep(D(A))$ for a factorizable Hopf algebra $A$.
\md
Then $\cK_1(\cD)=$ and $\cL_1(\cD)=?$.
\et
}
%%%%%%%%%
\subsection{On exact factorizations of factorizable Hopf algebras}
\bt\label{defcl}
Let $A$ be a modular (factorizable) Hopf algebra and $\Dc=\Rep(A//L)\subset \Rep(A)$ be a nondegenerate fusion subcategory. Suppose that $\Dc'=\Rep(A//C(L))$ for some normal left coideal subalgebra of $A$. Then $|C(L)|=\frac{|A|}{|L|}$, $L\cap C(L)=k$ and $L\otimes C(L)\ra A$ is an isomorphism of $A$-comodules.
\et

\bpf
The equality of dimensions follows applying Frobenius-Perron dimension formula, $\fp(\cD')\fp(\cD)=\fp(\cC)$.
Since $\Dc'\cap \Dc=\Vc$ it follows that $LC(L)=A$. On the other hand $\Rep(A)=\Dc \vee \Dc'\subset \Rep(A//(L\cap C(L)))$ and therefore $L\cap C(L)=k$.
\epf

We say that a normal coideal subalgebra $L$ of $A$ has a complement if there is another normal coideal subalgebra $C(L)$ such that such that the multiplication map $m:L \otimes C(L)\ra A$ is bijective.

\bc\label{chrnondeg}
Let $A$ be a modular Hopf algebra and $L$ be a normal coideal subalgebra of $A$. Then $\Rep(A//L)$ is nondegenerate if and only if $L$ has a normal coideal subalgebra complement $C(L)$. \ec

\bpf
If $\Rep(A//L)$ is nondegenerate then the previous lemma implies that $L$ has a complement $C(L)$. The converse of this follows from Theorem \ref{interms}.
\epf

\subsection{Some consequences of Theorem \ref{centraliz}}
\bp
Let $A$ be a \ssha and let $K(A)$ be its Hopf center. \thn 
\bne
\item 
$
(F^{-1}(\rep(A)_{ad}))'=<K(A)>
$
\item
$
\cd(\bar{G}(A))'=<V_{g, \eps}\;|\;g \in \bar{G}(A)>
$
\ene
\ep
\bpf
Note that $\cd(K(A))=F^{-1}(\rep(A)_{ad})$. Moreover $<\bar{G}(A)>=<V_{g, \eps}\;|\;g \in \bar{G}(A)>$.
\epf
\bc
Let $A$ be a \ssha with $C(A)$ commutative. \thn $ (F^{-1}(\rep(A)_{ad}))'=<V_{g, \eps}\;|\;g \in \bar{G}(A)>$.
\ec
\subsection {Some factorization results}
\bn{cor}
The fusion subcategory $\cD(K,\;L,\;\X,\;\psi)$ is non degenerate if and only if $\X$, $\psi$ are trivial and $A=K\otimes L$ as Hopf algebras.
\ec
\blue{Compare with Mueger and DNW}
\bpf
Dimension argument.
\epf
\bt\label{first factorization}
 If $A\cong K \otimes L$  as Hopf algebras then $D(A)\cong B(K,\;L)\otimes B(L,\;K)$ as Hopf algebras. %In particular $\Dc(K,\;L)'=\Dc(L,\;K)$.
 \et
\bpf
Dimension computations.
\epf
\subsection{The fusion subcategory generated by all normal coideals}
\blue{As consequence of equality in kernels}
\bc\label{main}
The fusion subcategory of $\Rep(D(A))$ generated by all normal coideals of $A$ coincides with $\Rep(D(A)//K(A))$.
\ec
\bn{proof}
By Theorem \ref{leftkr} one has that $ \LKER_A(A_{ad})= K(A)$. Then one can apply Brauer's generalized theorem \ref{charofim}.
\end{proof}
Denote by $\bar{G}(A)$ the set of all central grouplike elements of $A$. 

\bc Let $A$ be a semisimple Hopf algebra with the character ring $C(A)$ commutative (e.g. $A$ quasi-triangular). Then the fusion subcategory of normal coideal contains $\Rep \;D(A)_{\ad}$.
\ec

\bpf
If $C(A)$ commutative then $K(A)=k\bar{G}(A)$, see \cite{GN}. Thus one has $K(A)\subset K(D(A))=k\bar{G}(D(A))$. This implies  $\rep(D(A)//K(A))\supseteq \rep(D(A)//K(D(A)))=\Rep \;D(A)_{\ad}$.
\epf

\bn{cor}
Let $K$ be a normal Hopf subalgebra of $A$. Then 
$$\HKer_{D(A)}(K)=(A//K)^*\bowtie \HKer_A(K).$$ 
\end{cor}
\bpf
It is easy to see that the largest subcoalgebra of $D(A)$ contained in $\lker_{D(A)}(K)$ coincides to  $(A//K)^*\bowtie \HKer_A(K).$\epf

\subsection{On the fusion subcategory $\cd(A)$. Modules with trivial $A$-action.}

\bp\label{trivactcentr}
Let $V$ be a $D(A)$-module receiving a trivial $A$-action. Then $V$ centralizes another $D(A)$-module $M$ if and only if the subcoalgebra $C_V$ associated to $V$ as an $A$-comodule satisfies $C_V \subseteq \HKer_A(M)$.

In particular any two $D(A)$-modules receiving trivial $A$-action centralize each other.
\ep

\bpf
 Let $V$ be an $A$-module receiving trivial $A$-action. Then one has that $c_{V,\;M}:V\ot M\ra M\ot V$ is given by $v \ot a\mapsto a\ot v$. Thus $c_{V,\;M}c_{M,\;V}$ is given by $v\ot m\mapsto \sum_{i}e_i^*v \ot e_im=v_0\ot v_1m$.

If both $V$ and $M$ have trivial action then both braidings $c_{V,\;M}$ and $c_{M,\;V}$ are the twist map and therefore their composition is the identity.
\epf

\bc\lb{centralizer}
Let $V$ be an $A$-module receiving trivial $A$-action. Then $$<V>'=\Rep(D(A)//L(C_V)).$$
\ec
\bpf
\epf

\blue{It follows from the previous proposition that $\cd(A)$ is a symmetric category. Following \cite{ENO22} one has that $\cd(A)=\rep(D(A)//K(A))$. Thus $\cd(A)=<A>.$}
\subsubsection{The result from \cite{ENO2}} Suppose that $A$ is a semisimple Hopf algebra with $K(A)=kG^*$ for some group $G$. Then following Proposition 2.9 of \cite{ENO2} one has that $\Rep(G)\subset \Rep(D(A))$. If $M$ is an irreducible representation of $G$ then $M$ becomes a $D(A)$-module with the trivial $A$-action and the $A^*$-action given by $f.m=\sum_{g \in G}f(p_g)g.m$.

Moreover it is shown in the same paper that $$\Rep(G)'=\Rep(D(A)//K(A)).$$

It follows that $\Rep(G)=\Rep(D(A)//K(A))^*=<K(A)>$.
\red{\bn{rem} There is fusion subcategory of $D(H_8)$ which is not
normal. The one generated by a simple one dimensional simple which
is not central.{\bf other words, cite your paper}
\end{rem} }
\bp
Suppose that $L$ is a normal left coideal of $A$ and $M:=\lker_A(L)$. Then $L$ and $M$ centralize each other when regarded as $D(A)$-modules via Equation \ref{zhu} if and only if
\beq
Sm_4l_3m_3\ot Sm_2Sl_2m_5l_1m_1=l \ot m
\eeq
\ep

\bpf

One has the following formula for $c_{M,L}c_{L,M}$
\beq
c_{M,L}c_{L,M}(l \ot m)=Sm_4l_3m_3\ot Sm_2Sl_2m_5l_1m_1
\eeq

By symmetry one  has:
\beq
c_{L,M}c_{M,L}(m \ot l)=Sl_4m_3l_3\ot Sl_2Sm_2l_5m_1l_1
\eeq
\epf
\blue{Clearly $\lker_A(M)\supset L$. Do we have equality here?}

\blue{The dual module is given by $A^*$ with the adjoint action of $A^*$ and the action of $A$ given by $[a.f](x)=f(Sax)$.}

\blue {A different approach using the centralizer exercise from Kassel's book that $V_g \ot W_h$ is sent to $W_h\ot V_{hgh^{-1}}$. Thus if they centralize then the supports  commute elementwise.} 
%%%%%%%%
\bc
If $L$ and $M$ are two normal coalgebras then they centralize each other. if and only if one commutes with the dual of the other.
\ec
\subsection{Few more results on the centralizer}
\bt\lb{Mueger's centrliazer}
One has that $$\cD(<L\psi(x)>, \;\cap_{x\in \cX}\ker_K(x))\subseteq \cD( K, L,\;\cX,\psi)'.$$
\et
\bb{proof of this}
\bn{defn}
Let $K$ and $L$ be two Hopf subalgebras of $A$ and define $$D(K,\;L):=(A//K)_l^*\bwt L.$$ Let $\cD(K,\;L)$ be the full abelian subcategory of $\Rep(D(A))$ whose objects receive trivial action from $D(K,\;L)$.
\end{defn}

\bpf
The above inclusion is equivalent to $\cD( K, \;L,\;\cX,\psi)\subseteq \cD(\cap_{x\in \cX}\ker_K(x),\;<L\psi(x)>)$.
\epf

\bp
If $\cD(M,\;N)\subset \cD(K,\; L,\;\cX,\psi)'$ then \\$M
\supseteq <L\psi(x)>$ and $N\subseteq \cap_{x\in \cX}\ker_K(x)$.
\ep
See what is the maximum $D(M,\;N)$ that is included there.

\subsection{Lagrangian subcategories}
A fusi�n subcategory $\cD$ of $\cC$ is called Lagrangian if $\cD'=\cD$.

\bt
If $\Rep(D(A))$ has a normal Lagrangian fusion subcategory then $A$ is a Kac algebra.
\et
\bpf
Suppose that $\cD(K,\;L,\cX,\;\psi)$ is Lagrangian. Then $\cD(K,\;L,\cX,\;\psi)'=\cD(K,\;L,\cX,\;\psi)$. But we have seen that $\cD(K,\;L,\cX,\;\psi)'=\cD(M,\;N,\cX,\;\psi')$ with $M \supseteq L$ and $N\subseteq K$. If $M=K$ and $N=L$ one has that   $K=L$. Since $D(K,\;K)$ is normal it follows that $K$ is commutative and $(A//K)^*$ is also commutative. Thus $A//K$ is a cocommutative algebra and $A$ is a Kac algebra.
\epf

\subsection{Examples of centralizing fusion subcategories for quantum doubles}
 \bt\lb{centrlizeronediml}
 Suppose that $V$ is a one dimensional representation of $D(A)$ with character $\psi=(\phi, g)$. Then $M$ centralizes $V$ if and only if $\phi\bwt g$ acts trivially on $M$.
 \et

 \bpf
It is easy to compute that $C_{M,\;V}:M\otimes V \ra V \otimes M$ is given by $m\ot v\ra (\phi \bwt g)m\ot v$.
 \epf
 
From Theorem \ref{fact} it follows that $\eta\bw g$ is a central group like element of $D(A)$. Thus the elements centralizing $V$ form a normal fusion subcategory $D_{(\eta,g)}$. It follows from Remark \ref{verlinde} that $<V>'$ is normal but then $<V>=<V>''$ is not normal.
%%%%%%%%%%%%% 

\bt
Suppose that $V$ is a one dimensional module of $D(A)$. Then
$V=(\eta,\;g)$ by theorem \ref{fact} and 
$$\HKer_{D(A)}(V)=(A//L(g))^*\bwt (A/L(\eta))^*+ other terms$$

where $L(g)=<h_1gSh_2>$. Thus if $g$ central then $L(g)=k<g>$. Then $g$ acts on $M$ by a root of unity and $\eta$ by the inverse of that root. Thus the center of $V$ is larger than the component wise span.
\et 

What about the converse?
\blue{Double centralizer does not work; counterexample from groups.}
%\blue{Compute first for the $D(G)$}

\subsection{Here}

For $K$ a Hopf subalgebra of $A$ denote by $\cD(K)$ to be the fusion subcategory of $\mtr{Rep}(D(A))$ generated by all modules having trivial $K$-action. Clearly $\cD(K)=\Rep(D(A)//L(K))$ where $L(K)$ is the smallest normal left coideal subalgebra of $D(A)$ containing $K$.
\section{Centralizers for fusion subcategories of Deligne products and their centralizers}
\subsection{Fusion subcategories of direct products of categories} First we need a preliminary result for the structure of fusion subcategories of a Deligne product of two fusion categories.\md 
Let $\mtc{C}^1$ and $\mtc{C}^2$ be two fusion categories. Identify
them with the corresponding fusion subcategories of $\mtc{C}^1\boxtimes\mtc{C}^2$ given by $\mtc{C}^1\boxtimes 1$ and $1 \boxtimes \mtc{C}^2$ respectively. Then every simple object of
$\mtc{C}^1\boxtimes\mtc{C}^2$ is of the form $X_1\boxtimes X_2$
where $X_i$ is a simple object of $\mtc{C}^i$.
%%%
Let $\mtc{D} \subset \mtc{C}^1\boxtimes\mtc{C}^2$ be a fusion
subcategory. Define $\mtc{L}^i(\mtc{D}):=\mtc{D}\cap \mtc{C}^i$,
$i=1,2$. Let also $\mtc{K}^1(\mtc{D})$ be the fusion subcategory
generated by all simple objects $X_1$ of $\mtc{C}^1$ such that $X_1
\boxtimes X_2 \in \mtc{D}$ for some simple object $X_2$ of
$\mtc{C}^2$. Similarly define the fusion subcategory
$\mtc{K}^2(\mtc{D})$. Clearly
$\mtc{L}^1(\mtc{D})\boxtimes\mtc{L}^2(\mtc{D})\subset \mtc{D}
\subset\mtc{K}^1(\mtc{D})\boxtimes\mtc{K}^1(\mtc{D})$.
The following Lemma is taken from \cite{dgno}.
\bt\label{dg}
Let $\mtc{D} \subset \mtc{C}^1\boxtimes\mtc{C}^2$ be a fusion subcategory. Then there is
a group $\mtc{X}$ and faithful $\mtc{X}$-gradings $\mtc{K}_i(\mtc{D})=\oplus_{x
\in \mtc{X}}\mtc{K}^i(\mtc{D})_x$ with trivial components
$\mtc{L}^i(\mtc{D})$, $i=1, 2$ such that
\begin{equation}\label{81}
    \mtc{D}=\oplus_{x \in \mtc{X}}\mtc{K}^1(\mtc{D})_x\boxtimes\mtc{K}^2(\mtc{D})_x.
\end{equation}
\et
\bt\lb{interchange}
Let $\cD$ be a fusion subcategory of a Deligne product of two modular categories $\cC_1$ and $\cC_2$. Then 
$\cK_i(\cD')=\cL_i(\cD)'$ and $\cL_i(\cD')=\cK_i(\cD)'.$
Moreover the group $\cX$ is the same for $\cd$ and $\cd'$ and there is the following relation between gradings:
\et

\bp
If $\cD$ is Lagrangian fusion subcategory of $\cc_{1}\bx  \cc_{2}$ then $K_i(\cD)'=L_i(\cD)$ and $L_i(\cD)'=K_i(\cD).$\ep

\bp Let $\cc_{1}$ and $\cc_{2}$ be two braided fusion categories and form $\cc_{1}\bx \cc_{2}$.
Let $\cD$ and $\cE$ be two fusion subcategories of $\cC_1\bx \cC_2$. Then one has the following equalities

1) $\cK_i(\cD\cap \cE)=\cK_i(\cD)\cap \cK_i(\E)$

2) $\cK_i(\cD\vee \cE)=\cK_i(\cD)\vee \cK_i(\E)$

3) $\cL_i(\cD\cap \cE)=\cL_i(\cD)\cap \cL_i(\E)$

4) $\cL_i(\cD\vee \cE)=\cL_i(\cD)\vee \cL_i(\E)$
\ep

\bpf
The proof of the first two items is straightforward. For the last two items use Theorem \ref{interchange}.
Indeed one has that
\begin{eqnarray*}
 \cL_i(\cD\cap \cE)'&=&\cK_i((\cD\cap \cE)')=\cK_i(\cD' \vee \cE')=\\&=&\cK_i(\cD' ) \vee \cK_i(\cE')=\\&=&\cL_i(\cD)'\vee \cL_i(\cE)'.
\end{eqnarray*}
 Thus $\cL_i(\cD\cap \cE)=\cL_i(\cD)\cap \cL_i(\E).$
Similarly it can be proven that $\cL_i(\cD\vee \cE)=\cL_i(\cD)\vee \cL_i(\E).$
\epf

It is easy to see the following fact which also follows from a more general result for modular categories, namely that $\cZ(\cC)\cong \cC\boxtimes \cC^{\mtr{rev}}$ as braided fusion categories.

% now I have a formula for the M\"uger's centralizer

\subsection{The centralizer of a normal fusion subcategory of a modular category}
\bt\lb{form for normal}
Let $A$ be a factorizable semisimple Hopf algebra and $K$ be a normal Hopf subalgebra of $A$. Then
\beq\lb{Muegers}
\Rep(A//K)'\subseteq \cK_1(\Rep(D(A)//(A//K)^*\bwt \LKER_A(K)))
\eeq
\et

\bpf
As above let $\cD(K)$ be the fusion subcategory of $\Rep(D(A))$ with the property that $K$ acts trivially on each of its objects. Using Proposition \ref{justonecomp} it follows that 
\beqn
\cd(K)'\subseteq \rep(D(A)//(A//K)^*\bwt \LKER_A(K))
\eeqn
Note that by Proposition \ref{restrfact} one has $\cL_1(\cd(K))=\rep(A//K)$.
Thus Theorem \ref{interchange} gives that
\beqn\rep(A//K)' \subseteq \cK_1((\Rep(D(A)//(A//K)^*\bwt \LKER_A(K))))\eeqn 
and the proof is complete.\epf

Let $K$ be a normal Hopf subalgebra of $A$ and $K_{st}:=\cap_{x \in \mtr{St}_A(K)}\lker_K(x)$.

\bp
Let $A$ be a factorizable semisimple Hopf algebra and $K$ be a normal Hopf subalgebra of $A$. Then
\beq 
\Rep(A//L(K_{st}))'\subseteq \cK_1(\Rep(D(A)//(A//K)^*\bwt \LKER_A(K)))
\eeq
\ep
\bpf
Using Theorem \ref{interchange} \beqn\cK_1(\cD(K))'=\cL_1(\Rep(D(A)//(A//K)^*\bwt \LKER_A(K))). \eeqn But using Proposition \ref{restrfact} it is easy to see that $\cK_1(\cD(K))=\Rep(A//L(K_{st}))$.
\epf
\blue{\bt\lb{genchr}
Let $\cD:=\cD(K,\;L,\;\cX,\;\psi)$ be a a normal fusion subcategory of $\Rep(D(A))$ for a factorizable Hopf algebra $A$.
\md
Then $\cK_1(\cD)=$ and $\cL_1(\cD)=?$.
\et
}
%%%%%%%%%
\subsection{On exact factorizations of factorizable Hopf algebras}
\bt\label{defcl}
Let $A$ be a modular (factorizable) Hopf algebra and $\Dc=\Rep(A//L)\subset \Rep(A)$ be a nondegenerate fusion subcategory. Suppose that $\Dc'=\Rep(A//C(L))$ for some normal left coideal subalgebra of $A$. Then $|C(L)|=\frac{|A|}{|L|}$, $L\cap C(L)=k$ and $L\otimes C(L)\ra A$ is an isomorphism of $A$-comodules.
\et

\bpf
The equality of dimensions follows applying Frobenius-Perron dimension formula, $\fp(\cD')\fp(\cD)=\fp(\cC)$.
Since $\Dc'\cap \Dc=\Vc$ it follows that $LC(L)=A$. On the other hand $\Rep(A)=\Dc \vee \Dc'\subset \Rep(A//(L\cap C(L)))$ and therefore $L\cap C(L)=k$.
\epf

We say that a normal coideal subalgebra $L$ of $A$ has a complement if there is another normal coideal subalgebra $C(L)$ such that such that the multiplication map $m:L \otimes C(L)\ra A$ is bijective.

\bc\label{chrnondeg}
Let $A$ be a modular Hopf algebra and $L$ be a normal coideal subalgebra of $A$. Then $\Rep(A//L)$ is nondegenerate if and only if $L$ has a normal coideal subalgebra complement $C(L)$. \ec

\bpf
If $\Rep(A//L)$ is nondegenerate then the previous lemma implies that $L$ has a complement $C(L)$. The converse of this follows from Theorem \ref{interms}.
\epf

\subsection{Some consequences of Theorem \ref{centraliz}}
\bp
Let $A$ be a \ssha and let $K(A)$ be its Hopf center. \thn 
\bne
\item 
$
(F^{-1}(\rep(A)_{ad}))'=<K(A)>
$
\item
$
\cd(\bar{G}(A))'=<V_{g, \eps}\;|\;g \in \bar{G}(A)>
$
\ene
\ep
\bpf
Note that $\cd(K(A))=F^{-1}(\rep(A)_{ad})$. Moreover $<\bar{G}(A)>=<V_{g, \eps}\;|\;g \in \bar{G}(A)>$.
\epf
\bc
Let $A$ be a \ssha with $C(A)$ commutative. \thn $ (F^{-1}(\rep(A)_{ad}))'=<V_{g, \eps}\;|\;g \in \bar{G}(A)>$.
\ec
\subsection {Some factorization results}
\bn{cor}
The fusion subcategory $\cD(K,\;L,\;\X,\;\psi)$ is non degenerate if and only if $\X$, $\psi$ are trivial and $A=K\otimes L$ as Hopf algebras.
\ec
\blue{Compare with Mueger and DNW}
\bpf
Dimension argument.
\epf
\bt\label{first factorization}
 If $A\cong K \otimes L$  as Hopf algebras then $D(A)\cong B(K,\;L)\otimes B(L,\;K)$ as Hopf algebras. %In particular $\Dc(K,\;L)'=\Dc(L,\;K)$.
 \et
\bpf
Dimension computations.
\epf
\subsection{The fusion subcategory generated by all normal coideals}
\blue{As consequence of equality in kernels}
\bc\label{main}
The fusion subcategory of $\Rep(D(A))$ generated by all normal coideals of $A$ coincides with $\Rep(D(A)//K(A))$.
\ec
\bn{proof}
By Theorem \ref{leftkr} one has that $ \LKER_A(A_{ad})= K(A)$. Then one can apply Brauer's generalized theorem \ref{charofim}.
\end{proof}
Denote by $\bar{G}(A)$ the set of all central grouplike elements of $A$. 

\bc Let $A$ be a semisimple Hopf algebra with the character ring $C(A)$ commutative (e.g. $A$ quasi-triangular). Then the fusion subcategory of normal coideal contains $\Rep \;D(A)_{\ad}$.
\ec

\bpf
If $C(A)$ commutative then $K(A)=k\bar{G}(A)$, see \cite{GN}. Thus one has $K(A)\subset K(D(A))=k\bar{G}(D(A))$. This implies  $\rep(D(A)//K(A))\supseteq \rep(D(A)//K(D(A)))=\Rep \;D(A)_{\ad}$.
\epf

\bn{cor}
Let $K$ be a normal Hopf subalgebra of $A$. Then 
$$\HKer_{D(A)}(K)=(A//K)^*\bowtie \HKer_A(K).$$ 
\end{cor}
\bpf
It is easy to see that the largest subcoalgebra of $D(A)$ contained in $\lker_{D(A)}(K)$ coincides to  $(A//K)^*\bowtie \HKer_A(K).$\epf

\subsection{On the fusion subcategory $\cd(A)$. Modules with trivial $A$-action.}

\bp\label{trivactcentr}
Let $V$ be a $D(A)$-module receiving a trivial $A$-action. Then $V$ centralizes another $D(A)$-module $M$ if and only if the subcoalgebra $C_V$ associated to $V$ as an $A$-comodule satisfies $C_V \subseteq \HKer_A(M)$.

In particular any two $D(A)$-modules receiving trivial $A$-action centralize each other.
\ep

\bpf
 Let $V$ be an $A$-module receiving trivial $A$-action. Then one has that $c_{V,\;M}:V\ot M\ra M\ot V$ is given by $v \ot a\mapsto a\ot v$. Thus $c_{V,\;M}c_{M,\;V}$ is given by $v\ot m\mapsto \sum_{i}e_i^*v \ot e_im=v_0\ot v_1m$.

If both $V$ and $M$ have trivial action then both braidings $c_{V,\;M}$ and $c_{M,\;V}$ are the twist map and therefore their composition is the identity.
\epf

\bc\lb{centralizer}
Let $V$ be an $A$-module receiving trivial $A$-action. Then $$<V>'=\Rep(D(A)//L(C_V)).$$
\ec
\bpf
\epf

\blue{It follows from the previous proposition that $\cd(A)$ is a symmetric category. Following \cite{ENO22} one has that $\cd(A)=\rep(D(A)//K(A))$. Thus $\cd(A)=<A>.$}
\subsubsection{The result from \cite{ENO2}} Suppose that $A$ is a semisimple Hopf algebra with $K(A)=kG^*$ for some group $G$. Then following Proposition 2.9 of \cite{ENO2} one has that $\Rep(G)\subset \Rep(D(A))$. If $M$ is an irreducible representation of $G$ then $M$ becomes a $D(A)$-module with the trivial $A$-action and the $A^*$-action given by $f.m=\sum_{g \in G}f(p_g)g.m$.

Moreover it is shown in the same paper that $$\Rep(G)'=\Rep(D(A)//K(A)).$$

It follows that $\Rep(G)=\Rep(D(A)//K(A))^*=<K(A)>$.
\red{\bn{rem} There is fusion subcategory of $D(H_8)$ which is not
normal. The one generated by a simple one dimensional simple which
is not central.{\bf other words, cite your paper}
\end{rem} }
\bp
Suppose that $L$ is a normal left coideal of $A$ and $M:=\lker_A(L)$. Then $L$ and $M$ centralize each other when regarded as $D(A)$-modules via Equation \ref{zhu} if and only if
\beq
Sm_4l_3m_3\ot Sm_2Sl_2m_5l_1m_1=l \ot m
\eeq
\ep

\bpf

One has the following formula for $c_{M,L}c_{L,M}$
\beq
c_{M,L}c_{L,M}(l \ot m)=Sm_4l_3m_3\ot Sm_2Sl_2m_5l_1m_1
\eeq

By symmetry one  has:
\beq
c_{L,M}c_{M,L}(m \ot l)=Sl_4m_3l_3\ot Sl_2Sm_2l_5m_1l_1
\eeq
\epf
\blue{Clearly $\lker_A(M)\supset L$. Do we have equality here?}

\blue{The dual module is given by $A^*$ with the adjoint action of $A^*$ and the action of $A$ given by $[a.f](x)=f(Sax)$.}

\blue {A different approach using the centralizer exercise from Kassel's book that $V_g \ot W_h$ is sent to $W_h\ot V_{hgh^{-1}}$. Thus if they centralize then the supports  commute elementwise.} 
%%%%%%%%
\bc
If $L$ and $M$ are two normal coalgebras then they centralize each other. if and only if one commutes with the dual of the other.
\ec
\subsection{Few more results on the centralizer}
\bt\lb{Mueger's centrliazer}
One has that $$\cD(<L\psi(x)>, \;\cap_{x\in \cX}\ker_K(x))\subseteq \cD( K, L,\;\cX,\psi)'.$$
\et
\bb{proof of this}
\bn{defn}
Let $K$ and $L$ be two Hopf subalgebras of $A$ and define $$D(K,\;L):=(A//K)_l^*\bwt L.$$ Let $\cD(K,\;L)$ be the full abelian subcategory of $\Rep(D(A))$ whose objects receive trivial action from $D(K,\;L)$.
\end{defn}

\bpf
The above inclusion is equivalent to $\cD( K, \;L,\;\cX,\psi)\subseteq \cD(\cap_{x\in \cX}\ker_K(x),\;<L\psi(x)>)$.
\epf

\bp
If $\cD(M,\;N)\subset \cD(K,\; L,\;\cX,\psi)'$ then \\$M
\supseteq <L\psi(x)>$ and $N\subseteq \cap_{x\in \cX}\ker_K(x)$.
\ep
See what is the maximum $D(M,\;N)$ that is included there.

\subsection{Lagrangian subcategories}
A fusi�n subcategory $\cD$ of $\cC$ is called Lagrangian if $\cD'=\cD$.

\bt
If $\Rep(D(A))$ has a normal Lagrangian fusion subcategory then $A$ is a Kac algebra.
\et
\bpf
Suppose that $\cD(K,\;L,\cX,\;\psi)$ is Lagrangian. Then $\cD(K,\;L,\cX,\;\psi)'=\cD(K,\;L,\cX,\;\psi)$. But we have seen that $\cD(K,\;L,\cX,\;\psi)'=\cD(M,\;N,\cX,\;\psi')$ with $M \supseteq L$ and $N\subseteq K$. If $M=K$ and $N=L$ one has that   $K=L$. Since $D(K,\;K)$ is normal it follows that $K$ is commutative and $(A//K)^*$ is also commutative. Thus $A//K$ is a cocommutative algebra and $A$ is a Kac algebra.
\epf

\subsection{Examples of centralizing fusion subcategories for quantum doubles}
 \bt\lb{centrlizeronediml}
 Suppose that $V$ is a one dimensional representation of $D(A)$ with character $\psi=(\phi, g)$. Then $M$ centralizes $V$ if and only if $\phi\bwt g$ acts trivially on $M$.
 \et

 \bpf
It is easy to compute that $C_{M,\;V}:M\otimes V \ra V \otimes M$ is given by $m\ot v\ra (\phi \bwt g)m\ot v$.
 \epf
 
From Theorem \ref{fact} it follows that $\eta\bw g$ is a central group like element of $D(A)$. Thus the elements centralizing $V$ form a normal fusion subcategory $D_{(\eta,g)}$. It follows from Remark \ref{verlinde} that $<V>'$ is normal but then $<V>=<V>''$ is not normal.
%%%%%%%%%%%%% 

\bt
Suppose that $V$ is a one dimensional module of $D(A)$. Then
$V=(\eta,\;g)$ by theorem \ref{fact} and 
$$\HKer_{D(A)}(V)=(A//L(g))^*\bwt (A/L(\eta))^*+ other terms$$

where $L(g)=<h_1gSh_2>$. Thus if $g$ central then $L(g)=k<g>$. Then $g$ acts on $M$ by a root of unity and $\eta$ by the inverse of that root. Thus the center of $V$ is larger than the component wise span.
\et 

What about the converse?
\blue{Double centralizer does not work; counterexample from groups.}
%\blue{Compute first for the $D(G)$}

\subsection{Here}

For $K$ a Hopf subalgebra of $A$ denote by $\cD(K)$ to be the fusion subcategory of $\mtr{Rep}(D(A))$ generated by all modules having trivial $K$-action. Clearly $\cD(K)=\Rep(D(A)//L(K))$ where $L(K)$ is the smallest normal left coideal subalgebra of $D(A)$ containing $K$.
\section{Centralizers for fusion subcategories of Deligne products and their centralizers}
\subsection{Fusion subcategories of direct products of categories} First we need a preliminary result for the structure of fusion subcategories of a Deligne product of two fusion categories.\md 
Let $\mtc{C}^1$ and $\mtc{C}^2$ be two fusion categories. Identify
them with the corresponding fusion subcategories of $\mtc{C}^1\boxtimes\mtc{C}^2$ given by $\mtc{C}^1\boxtimes 1$ and $1 \boxtimes \mtc{C}^2$ respectively. Then every simple object of
$\mtc{C}^1\boxtimes\mtc{C}^2$ is of the form $X_1\boxtimes X_2$
where $X_i$ is a simple object of $\mtc{C}^i$.
%%%
Let $\mtc{D} \subset \mtc{C}^1\boxtimes\mtc{C}^2$ be a fusion
subcategory. Define $\mtc{L}^i(\mtc{D}):=\mtc{D}\cap \mtc{C}^i$,
$i=1,2$. Let also $\mtc{K}^1(\mtc{D})$ be the fusion subcategory
generated by all simple objects $X_1$ of $\mtc{C}^1$ such that $X_1
\boxtimes X_2 \in \mtc{D}$ for some simple object $X_2$ of
$\mtc{C}^2$. Similarly define the fusion subcategory
$\mtc{K}^2(\mtc{D})$. Clearly
$\mtc{L}^1(\mtc{D})\boxtimes\mtc{L}^2(\mtc{D})\subset \mtc{D}
\subset\mtc{K}^1(\mtc{D})\boxtimes\mtc{K}^1(\mtc{D})$.
The following Lemma is taken from \cite{dgno}.
\bt\label{dg}
Let $\mtc{D} \subset \mtc{C}^1\boxtimes\mtc{C}^2$ be a fusion subcategory. Then there is
a group $\mtc{X}$ and faithful $\mtc{X}$-gradings $\mtc{K}_i(\mtc{D})=\oplus_{x
\in \mtc{X}}\mtc{K}^i(\mtc{D})_x$ with trivial components
$\mtc{L}^i(\mtc{D})$, $i=1, 2$ such that
\begin{equation}\label{81}
    \mtc{D}=\oplus_{x \in \mtc{X}}\mtc{K}^1(\mtc{D})_x\boxtimes\mtc{K}^2(\mtc{D})_x.
\end{equation}
\et
\bt\lb{interchange}
Let $\cD$ be a fusion subcategory of a Deligne product of two modular categories $\cC_1$ and $\cC_2$. Then 
$\cK_i(\cD')=\cL_i(\cD)'$ and $\cL_i(\cD')=\cK_i(\cD)'.$
Moreover the group $\cX$ is the same for $\cd$ and $\cd'$ and there is the following relation between gradings:
\et

\bp
If $\cD$ is Lagrangian fusion subcategory of $\cc_{1}\bx  \cc_{2}$ then $K_i(\cD)'=L_i(\cD)$ and $L_i(\cD)'=K_i(\cD).$\ep

\bp Let $\cc_{1}$ and $\cc_{2}$ be two braided fusion categories and form $\cc_{1}\bx \cc_{2}$.
Let $\cD$ and $\cE$ be two fusion subcategories of $\cC_1\bx \cC_2$. Then one has the following equalities

1) $\cK_i(\cD\cap \cE)=\cK_i(\cD)\cap \cK_i(\E)$

2) $\cK_i(\cD\vee \cE)=\cK_i(\cD)\vee \cK_i(\E)$

3) $\cL_i(\cD\cap \cE)=\cL_i(\cD)\cap \cL_i(\E)$

4) $\cL_i(\cD\vee \cE)=\cL_i(\cD)\vee \cL_i(\E)$
\ep

\bpf
The proof of the first two items is straightforward. For the last two items use Theorem \ref{interchange}.
Indeed one has that
\begin{eqnarray*}
 \cL_i(\cD\cap \cE)'&=&\cK_i((\cD\cap \cE)')=\cK_i(\cD' \vee \cE')=\\&=&\cK_i(\cD' ) \vee \cK_i(\cE')=\\&=&\cL_i(\cD)'\vee \cL_i(\cE)'.
\end{eqnarray*}
 Thus $\cL_i(\cD\cap \cE)=\cL_i(\cD)\cap \cL_i(\E).$
Similarly it can be proven that $\cL_i(\cD\vee \cE)=\cL_i(\cD)\vee \cL_i(\E).$
\epf

It is easy to see the following fact which also follows from a more general result for modular categories, namely that $\cZ(\cC)\cong \cC\boxtimes \cC^{\mtr{rev}}$ as braided fusion categories.

% now I have a formula for the M\"uger's centralizer

\subsection{The centralizer of a normal fusion subcategory of a modular category}
\bt\lb{form for normal}
Let $A$ be a factorizable semisimple Hopf algebra and $K$ be a normal Hopf subalgebra of $A$. Then
\beq\lb{Muegers}
\Rep(A//K)'\subseteq \cK_1(\Rep(D(A)//(A//K)^*\bwt \LKER_A(K)))
\eeq
\et

\bpf
As above let $\cD(K)$ be the fusion subcategory of $\Rep(D(A))$ with the property that $K$ acts trivially on each of its objects. Using Proposition \ref{justonecomp} it follows that 
\beqn
\cd(K)'\subseteq \rep(D(A)//(A//K)^*\bwt \LKER_A(K))
\eeqn
Note that by Proposition \ref{restrfact} one has $\cL_1(\cd(K))=\rep(A//K)$.
Thus Theorem \ref{interchange} gives that
\beqn\rep(A//K)' \subseteq \cK_1((\Rep(D(A)//(A//K)^*\bwt \LKER_A(K))))\eeqn 
and the proof is complete.\epf

Let $K$ be a normal Hopf subalgebra of $A$ and $K_{st}:=\cap_{x \in \mtr{St}_A(K)}\lker_K(x)$.

\bp
Let $A$ be a factorizable semisimple Hopf algebra and $K$ be a normal Hopf subalgebra of $A$. Then
\beq 
\Rep(A//L(K_{st}))'\subseteq \cK_1(\Rep(D(A)//(A//K)^*\bwt \LKER_A(K)))
\eeq
\ep
\bpf
Using Theorem \ref{interchange} \beqn\cK_1(\cD(K))'=\cL_1(\Rep(D(A)//(A//K)^*\bwt \LKER_A(K))). \eeqn But using Proposition \ref{restrfact} it is easy to see that $\cK_1(\cD(K))=\Rep(A//L(K_{st}))$.
\epf
\blue{\bt\lb{genchr}
Let $\cD:=\cD(K,\;L,\;\cX,\;\psi)$ be a a normal fusion subcategory of $\Rep(D(A))$ for a factorizable Hopf algebra $A$.
\md
Then $\cK_1(\cD)=$ and $\cL_1(\cD)=?$.
\et
}
%%%%%%%%%
\subsection{On exact factorizations of factorizable Hopf algebras}
\bt\label{defcl}
Let $A$ be a modular (factorizable) Hopf algebra and $\Dc=\Rep(A//L)\subset \Rep(A)$ be a nondegenerate fusion subcategory. Suppose that $\Dc'=\Rep(A//C(L))$ for some normal left coideal subalgebra of $A$. Then $|C(L)|=\frac{|A|}{|L|}$, $L\cap C(L)=k$ and $L\otimes C(L)\ra A$ is an isomorphism of $A$-comodules.
\et

\bpf
The equality of dimensions follows applying Frobenius-Perron dimension formula, $\fp(\cD')\fp(\cD)=\fp(\cC)$.
Since $\Dc'\cap \Dc=\Vc$ it follows that $LC(L)=A$. On the other hand $\Rep(A)=\Dc \vee \Dc'\subset \Rep(A//(L\cap C(L)))$ and therefore $L\cap C(L)=k$.
\epf

We say that a normal coideal subalgebra $L$ of $A$ has a complement if there is another normal coideal subalgebra $C(L)$ such that such that the multiplication map $m:L \otimes C(L)\ra A$ is bijective.

\bc\label{chrnondeg}
Let $A$ be a modular Hopf algebra and $L$ be a normal coideal subalgebra of $A$. Then $\Rep(A//L)$ is nondegenerate if and only if $L$ has a normal coideal subalgebra complement $C(L)$. \ec

\bpf
If $\Rep(A//L)$ is nondegenerate then the previous lemma implies that $L$ has a complement $C(L)$. The converse of this follows from Theorem \ref{interms}.
\epf

\subsection{Some consequences of Theorem \ref{centraliz}}
\bp
Let $A$ be a \ssha and let $K(A)$ be its Hopf center. \thn 
\bne
\item 
$
(F^{-1}(\rep(A)_{ad}))'=<K(A)>
$
\item
$
\cd(\bar{G}(A))'=<V_{g, \eps}\;|\;g \in \bar{G}(A)>
$
\ene
\ep
\bpf
Note that $\cd(K(A))=F^{-1}(\rep(A)_{ad})$. Moreover $<\bar{G}(A)>=<V_{g, \eps}\;|\;g \in \bar{G}(A)>$.
\epf
\bc
Let $A$ be a \ssha with $C(A)$ commutative. \thn $ (F^{-1}(\rep(A)_{ad}))'=<V_{g, \eps}\;|\;g \in \bar{G}(A)>$.
\ec
\subsection {Some factorization results}
\bn{cor}
The fusion subcategory $\cD(K,\;L,\;\X,\;\psi)$ is non degenerate if and only if $\X$, $\psi$ are trivial and $A=K\otimes L$ as Hopf algebras.
\ec
\blue{Compare with Mueger and DNW}
\bpf
Dimension argument.
\epf
\bt\label{first factorization}
 If $A\cong K \otimes L$  as Hopf algebras then $D(A)\cong B(K,\;L)\otimes B(L,\;K)$ as Hopf algebras. %In particular $\Dc(K,\;L)'=\Dc(L,\;K)$.
 \et
\bpf
Dimension computations.
\epf
\subsection{The fusion subcategory generated by all normal coideals}
\blue{As consequence of equality in kernels}
\bc\label{main}
The fusion subcategory of $\Rep(D(A))$ generated by all normal coideals of $A$ coincides with $\Rep(D(A)//K(A))$.
\ec
\bn{proof}
By Theorem \ref{leftkr} one has that $ \LKER_A(A_{ad})= K(A)$. Then one can apply Brauer's generalized theorem \ref{charofim}.
\end{proof}
Denote by $\bar{G}(A)$ the set of all central grouplike elements of $A$. 

\bc Let $A$ be a semisimple Hopf algebra with the character ring $C(A)$ commutative (e.g. $A$ quasi-triangular). Then the fusion subcategory of normal coideal contains $\Rep \;D(A)_{\ad}$.
\ec

\bpf
If $C(A)$ commutative then $K(A)=k\bar{G}(A)$, see \cite{GN}. Thus one has $K(A)\subset K(D(A))=k\bar{G}(D(A))$. This implies  $\rep(D(A)//K(A))\supseteq \rep(D(A)//K(D(A)))=\Rep \;D(A)_{\ad}$.
\epf

\bn{cor}
Let $K$ be a normal Hopf subalgebra of $A$. Then 
$$\HKer_{D(A)}(K)=(A//K)^*\bowtie \HKer_A(K).$$ 
\end{cor}
\bpf
It is easy to see that the largest subcoalgebra of $D(A)$ contained in $\lker_{D(A)}(K)$ coincides to  $(A//K)^*\bowtie \HKer_A(K).$\epf

\subsection{On the fusion subcategory $\cd(A)$. Modules with trivial $A$-action.}

\bp\label{trivactcentr}
Let $V$ be a $D(A)$-module receiving a trivial $A$-action. Then $V$ centralizes another $D(A)$-module $M$ if and only if the subcoalgebra $C_V$ associated to $V$ as an $A$-comodule satisfies $C_V \subseteq \HKer_A(M)$.

In particular any two $D(A)$-modules receiving trivial $A$-action centralize each other.
\ep

\bpf
 Let $V$ be an $A$-module receiving trivial $A$-action. Then one has that $c_{V,\;M}:V\ot M\ra M\ot V$ is given by $v \ot a\mapsto a\ot v$. Thus $c_{V,\;M}c_{M,\;V}$ is given by $v\ot m\mapsto \sum_{i}e_i^*v \ot e_im=v_0\ot v_1m$.

If both $V$ and $M$ have trivial action then both braidings $c_{V,\;M}$ and $c_{M,\;V}$ are the twist map and therefore their composition is the identity.
\epf

\bc\lb{centralizer}
Let $V$ be an $A$-module receiving trivial $A$-action. Then $$<V>'=\Rep(D(A)//L(C_V)).$$
\ec
\bpf
\epf

\blue{It follows from the previous proposition that $\cd(A)$ is a symmetric category. Following \cite{ENO22} one has that $\cd(A)=\rep(D(A)//K(A))$. Thus $\cd(A)=<A>.$}
\subsubsection{The result from \cite{ENO2}} Suppose that $A$ is a semisimple Hopf algebra with $K(A)=kG^*$ for some group $G$. Then following Proposition 2.9 of \cite{ENO2} one has that $\Rep(G)\subset \Rep(D(A))$. If $M$ is an irreducible representation of $G$ then $M$ becomes a $D(A)$-module with the trivial $A$-action and the $A^*$-action given by $f.m=\sum_{g \in G}f(p_g)g.m$.

Moreover it is shown in the same paper that $$\Rep(G)'=\Rep(D(A)//K(A)).$$

It follows that $\Rep(G)=\Rep(D(A)//K(A))^*=<K(A)>$.
\red{\bn{rem} There is fusion subcategory of $D(H_8)$ which is not
normal. The one generated by a simple one dimensional simple which
is not central.{\bf other words, cite your paper}
\end{rem} }
\bp
Suppose that $L$ is a normal left coideal of $A$ and $M:=\lker_A(L)$. Then $L$ and $M$ centralize each other when regarded as $D(A)$-modules via Equation \ref{zhu} if and only if
\beq
Sm_4l_3m_3\ot Sm_2Sl_2m_5l_1m_1=l \ot m
\eeq
\ep

\bpf

One has the following formula for $c_{M,L}c_{L,M}$
\beq
c_{M,L}c_{L,M}(l \ot m)=Sm_4l_3m_3\ot Sm_2Sl_2m_5l_1m_1
\eeq

By symmetry one  has:
\beq
c_{L,M}c_{M,L}(m \ot l)=Sl_4m_3l_3\ot Sl_2Sm_2l_5m_1l_1
\eeq
\epf
\blue{Clearly $\lker_A(M)\supset L$. Do we have equality here?}

\blue{The dual module is given by $A^*$ with the adjoint action of $A^*$ and the action of $A$ given by $[a.f](x)=f(Sax)$.}

\blue {A different approach using the centralizer exercise from Kassel's book that $V_g \ot W_h$ is sent to $W_h\ot V_{hgh^{-1}}$. Thus if they centralize then the supports  commute elementwise.} 
%%%%%%%%
\bc
If $L$ and $M$ are two normal coalgebras then they centralize each other. if and only if one commutes with the dual of the other.
\ec
\subsection{Few more results on the centralizer}
\bt\lb{Mueger's centrliazer}
One has that $$\cD(<L\psi(x)>, \;\cap_{x\in \cX}\ker_K(x))\subseteq \cD( K, L,\;\cX,\psi)'.$$
\et
\bb{proof of this}
\bn{defn}
Let $K$ and $L$ be two Hopf subalgebras of $A$ and define $$D(K,\;L):=(A//K)_l^*\bwt L.$$ Let $\cD(K,\;L)$ be the full abelian subcategory of $\Rep(D(A))$ whose objects receive trivial action from $D(K,\;L)$.
\end{defn}

\bpf
The above inclusion is equivalent to $\cD( K, \;L,\;\cX,\psi)\subseteq \cD(\cap_{x\in \cX}\ker_K(x),\;<L\psi(x)>)$.
\epf

\bp
If $\cD(M,\;N)\subset \cD(K,\; L,\;\cX,\psi)'$ then \\$M
\supseteq <L\psi(x)>$ and $N\subseteq \cap_{x\in \cX}\ker_K(x)$.
\ep
See what is the maximum $D(M,\;N)$ that is included there.

\subsection{Lagrangian subcategories}
A fusi�n subcategory $\cD$ of $\cC$ is called Lagrangian if $\cD'=\cD$.

\bt
If $\Rep(D(A))$ has a normal Lagrangian fusion subcategory then $A$ is a Kac algebra.
\et
\bpf
Suppose that $\cD(K,\;L,\cX,\;\psi)$ is Lagrangian. Then $\cD(K,\;L,\cX,\;\psi)'=\cD(K,\;L,\cX,\;\psi)$. But we have seen that $\cD(K,\;L,\cX,\;\psi)'=\cD(M,\;N,\cX,\;\psi')$ with $M \supseteq L$ and $N\subseteq K$. If $M=K$ and $N=L$ one has that   $K=L$. Since $D(K,\;K)$ is normal it follows that $K$ is commutative and $(A//K)^*$ is also commutative. Thus $A//K$ is a cocommutative algebra and $A$ is a Kac algebra.
\epf

\subsection{Examples of centralizing fusion subcategories for quantum doubles}
 \bt\lb{centrlizeronediml}
 Suppose that $V$ is a one dimensional representation of $D(A)$ with character $\psi=(\phi, g)$. Then $M$ centralizes $V$ if and only if $\phi\bwt g$ acts trivially on $M$.
 \et

 \bpf
It is easy to compute that $C_{M,\;V}:M\otimes V \ra V \otimes M$ is given by $m\ot v\ra (\phi \bwt g)m\ot v$.
 \epf
 
From Theorem \ref{fact} it follows that $\eta\bw g$ is a central group like element of $D(A)$. Thus the elements centralizing $V$ form a normal fusion subcategory $D_{(\eta,g)}$. It follows from Remark \ref{verlinde} that $<V>'$ is normal but then $<V>=<V>''$ is not normal.
%%%%%%%%%%%%% 

\bt
Suppose that $V$ is a one dimensional module of $D(A)$. Then
$V=(\eta,\;g)$ by theorem \ref{fact} and 
$$\HKer_{D(A)}(V)=(A//L(g))^*\bwt (A/L(\eta))^*+ other terms$$

where $L(g)=<h_1gSh_2>$. Thus if $g$ central then $L(g)=k<g>$. Then $g$ acts on $M$ by a root of unity and $\eta$ by the inverse of that root. Thus the center of $V$ is larger than the component wise span.
\et 

What about the converse?
\blue{Double centralizer does not work; counterexample from groups.}
%\blue{Compute first for the $D(G)$}

\section{Examples and applications}\lb{appl}
\stat 18
\bn{conj}
If $S(g, M)$ centralize $S(h, N)$ then $\cC(g)$ commutes elemtwise with $\cc(h)$.
\end{conj}
If $S(g, M)$ centralize $S(h, N)$ then $gh=hg$. Look at the induced map of the braiding on the restriction to $K(A)$.
\md
for $H< U(A)$ define $\cc(H)$ as the set of all simple modules supported on elements of $H$. Is $\cc(H)'=\cc(K)$ for some $K$ with $[K, H]=1$.
\stat 19
\mdn the simple modules in the $D(A)$-module $A$ are of the form $A\ot_{L(g)}k$.
\subsection{On the commutator category}
\bp
Let $K$ be a Hopf subalgebra of $A$ and denote $\cd(K)$ the fusion subcategory of $Rep(A)$ whose objects receive trivial action from $A$. hen 
\beq
\cd(K)^{co}=\rep(A/[A, K])
\eeq
\ep
\bpf
First note that $\rep(A/[A, K])$ consists of those objects whose restriction to $K$ is a homogenous one dimensional module. If $(am-ma).v=0$ then $L_{m}$ is a morphism of $A$-modules and therefore a scalar of identity.
\md
It follows that $\rep(A/[A, K])\subset \cd(K)^{co}$. Conversely if $V\ot V^{*}$ is trivial as $K$-module it follows that
$V\dw_{K}=\dim_{\kk} V M$ where $M$ is a one dimensional module of $K$.
\epf

\blue{If $A$ is factorizable what is the correspondent $\ch_g$? Apply it to various Hopf algebras and find fusion subcategories of dimension the index of $g$.}

\subsection{On the left kernel and the Hopf kernel of the adjoint action}
%%%%%%%%%%
Let $A$ be a semisimple Hopf algebra and let $A_{ad}:=A$ be the adjoint left $A$-module given by: $x.a=x_1aS(x_2)$ for all $x, a \in A$.

\bp\lb{ladj}
Suppose $M\subseteq A$ is an $A$-submodule of $A_{\ad}$. Then the left kernel $\LKER_A(M)$ is the largest left coideal subalgebra $R$ of $A$ such that $S(R)$ commutes elementwise with $M$. If $S(M)=M$ then $R$ is the largest left coideal subalgebra of $A$ that commutes element wise with $M$.
\ep
\bpf
Suppose that $r\in  \LKER_A(M)$. Note that by Equation \eqref{L} one has $r_{1}\ot r_{2}mS(r_{3})=r \ot m$
for all $m \in M$.
This gives that $mS(r)=S(r_1)(r_2mS(r_3))=S(r)m.$ Thus $S(\lker_A(M))$ commutes elementwise to $M$. Conversely if $S(R)$ commutes elementwise with $M$ then for any $r \in R$ one has $r_1\ot r_2mS(r_3)=r_1\ot r_2S(r_3)m=r\ot m$ and therefore $R \subseteq \LKER_A(M)$.
\epf
\br\label{hadj}
Similarly it can be proven that $\hker_A(M)$ is the largest Hopf subalgebra of $A$ that commutes elementwise with $M$.
\er
\bc\label{hopadj}
Suppose $M$ is a submodule of $A_{\ad}$. Then the Hopf kernel $\HKer_A(M)$ is the largest Hopf subalgebra $K$ of $A$ that commutes element wise with $M$.
\ec
%\bpf Straightforward from the previous Proposition.\epf

\bp \lb{leftcoinkernadj}Let $A$ be a semisimple Hopf algebra and $A_{ad}:=A$ be the adjoint $A$-module given by: $x.a=x_1aS(x_2)$.
Then 
\bne
\item
$\mtr{LKer}_{ _A}(A_{ad})$ is the largest central coideal subalgebra of $A$.
\item $\mtr{HKer}_{A_{ad}}$ is the largest Hopf central Hopf subalgebra of $A$.
\ene
 Moreover $\mtr{LKer}_{ _A}(A_{ad}) =\hker_{A}(A_{ad})$, the largest central Hopf subalgebra of $A$.
\ep

\bpf Apply previous proposition for $M=A_{ad}$. Moreover by \cite[Theorem 2.8]{GN}. one has that $\rep(A//K(A))=<A_{ad}>$ which shows that $\mtr{LKer}_{ _A}(A_{ad}) =\hker_{A}(A_{ad})$.\epf

%\subsection{The fusion subcategory generated by all left normal coideals}
\subsection{On a $D(A)$-module strcuture on $A$.}\label{mod}
Let $A$ be a semisimple Hopf algebra and $D(A)$ be its Drinfeld double. Then $A$ is a left $D(A)$-module via
\beq \lb{zhu}(f \bwt a).b=(a_1bS(a_2))\lh S^{-1}f\eeq
for all $f \in A^*$ and $a, b\in A$.
\md
The simple $D(A)$-submodules of $A$ are the minimal normal left coideals of $A$.  Also note that any normal left coideal subalgebra $L$ of $A$ can be regarded as a $D(A)$-submodule of $A$.

 \section{Appendix: Left kernels and normal left coideal subalgebras of a Drinfeld double}
In this section we give a theorem structure for coideal subalgebras of Drinfeld doubles.
\lb{leftkernelsinda}
\subsection{On the integrals of coideal subalgebras of a semisimple Hopf algebra}
\subsubsection{Definition of the invariant element }\cite{gmj} Let$L$ be a left coideal subalgebra of a semisimple Hopf algebra $A$. An element $x \in S$ is called invariant if $ls=\eps(l)s$ for all $s \in S$. It is shown in \cite{gmj} that the invariant elements is always unique up to a scalar. Moreover \cite[Theorem 2.1]{kacoid} implies that 
\beq\lb{invr}
L=x_{L}\lh A^{*}
\eeq
for any invariant element $x_{L}$. Moreover if $L$ is a normal coideal subalgebra then $x_{L}$ is central in $A$.
\bp\label{closetofrob}
Let $S$ be a coideal subalgebra of a Hopf algebra $A$. Then $$\rho(x_S)=\sum_{i=1}^sa_i\ot s_i$$ where $s_i$ is a basis of $S$ and $a_i$ are linearly independent.
\ep

\bpf
\blu{Use the fact that $S=A^*\rh x_S$ which was proved in Theorem 2.1 of \cite{kacoid}.}
\epf

\subsection{The invariant element of a coideal subalgebra}\label{integral}
Let $L$ be any left coideal subalgebra of a finite dimensional Hopf algebra $A$. Then $A$ is free over $L$ \cite{Sk} both as left and right $L$-module. Let $A=L\oplus(\oplus_{i=2}^rLx_i)$ be a decomposition of $A$ as a free left $L$-modules of rank one. Then the idempotent integral $\Lam$ of $A$ admits a decomposition: $$\Lam=y_S+\sum_{i=2}^ry_i$$ with $y_S \in S$ and $y_i \in Sx_i$ for $i \geq 2$. Then equation $s\Lam=\eps(s)\Lam$ implies that $sy_S=\epsilon(s)y_S$ for all $s \in S$. 
\subsection{Normal left coideal subalgebras of  $D(A)$} Since $A$ si also cosemisimple there is a decomposition $A=\bigoplus_{d \in \Irr(A^{*})}C_{d}$ where $C_{d}$ is the coefficient subcoalgebra of $A$ associated to the irreducible character $d \in \irr(A^{*})$.
Let $x^d_{ij}$ be a fixed comatrix basis for the subcoalgebra $C_d$. Recall that $\D(x^d_{ij})=\sum_k x^d_{ik}\otimes x^d_{kj}$ for all $i,j$.
Note that
\beq
a^{d}_{ij}\lh A^{d}_{uv}=\delta_{i, u}a^{d}_{vj}
\eeq
\bpf   Since $L$ is a left coideal it follows that $l\lh D(A)^{*} \in L$. One may obviously identify $D(A)^{*}\cong A\ot A^{*}$ as vector spaces. Let $A^{d}_{ij}$ be the dual basis of the matrix coefficient bases $\{ a_{ij}^{d}\}$. Therefore  $A\rh F_{dij}\ot a^d_{ij} =l \lh (A\ot A^{d}_{ii})\in L$ for any $d,i,j,k$. In particular $F_{dij}\ot x^d_{ij}\in L$ for any $d,i,j$.
\epf

\bp \label{descr}
Let L be a left coideal subalgebra of $D(A)$ and suppose that there are vector subspaces $B^d_{ij}$ of $A^*$ such that 
\beq\lb{eqs}
L=\bigoplus_{ijd}B^d_{ij}\bwt ka^d_{ij}
\eeq
Then $B^d_{ij}=B^d_{kj}$ for all $1\leq k \leq \eps(d)$. Moreover $L$ can be written as
\beq
L=\bigoplus_{i=1}^{s}(M_{d,i}\bwt L_{d,i})
\eeq
where $L_{d,i}$ are the simple left coideals in $A$ and $M_{d,i}$ are some right coideals in $A^*$. 
\ep
\blue{For coideal subalgebras of Kac algebras this is true from the arguments from \cite{coidsuabalg}}
\bpf Clearly there 
are vector subspaces $B^d_{ij}\subset A^*$ such that 
$
L=\oplus_{ijd}B^d_{ij}\bwt ka^d_{ij}.
$ Let $b^{d}_{ij}\in B^{d}_{ij}$ arbitrary chosen. Then it follows that 
\beq
b^{d}_{ij}\bwt a^{d}_{ij} \lh (A\ot A^{d}_{ik})=(A\rh b^{d}_{ij})
\bwt (a^{d}_{ij}\lh A^{d}_{ik})=(A\rh b^{d}_{ij})
\bwt a^{d}_{kj}
\eeq
This shows that $(A\rh b^{d}_{ij})\subseteq B_{vj}^{d}$ and in particular $B^{d}_{ij}\subseteq B^{d}_{kj}$. By symmetry one has that $B_{kj}^{d}\subseteq B^{d}_{ij}$ which shows that $B_{kj}^{d}=B^{d}_{ij}$.
\epf

It is easy now to check the following corollary:
\bc\lb{onesd} Let $L$ be any left coideal of a semisimple Drinfeld double $D(A)$. Suppose that $l=\sum_{i=1}^{s}f_{ij}^{d}\bwt a_{ij }^{d}\in L$ where $a^{d}_{ij}$ is the matrix coefficient $\kk$-linear basis of $A$. Then $(A\rh f_{ij }^{d})\ot (a_{ij}^{d}\lh A^{*}) \subset L$ for all $1\leq i, j\leq \eps(d)$ and all $d \in \irr(A^{*})$.
\ec
Let $x_L=\sum_{ijd}f^d_{ij}\bwt a^d_{ij}$ be the invariant element of $L$. It follows that $x_L\lh D(A)^*=L$. Identify $D(A)^*$ with $A \otimes A^*$ and let as above $A^d_{ij}$ be the dual basis of the comatrix basis $a^d_{ij}$.
\bn{conj}Conjecture: One has $B^{d}_{ij}= A\rh f^{d}_{ii}$ for all ${d,i,j}$.
\end{conj}
\bl For normal coideal subalgebras one has
\beq
(A\bwt B)\cap (A' \bwt B')=A\cap A' \bwt B\cap B'
\eeq
\el
\bpf
Denote $L:=(A\bwt B)\cap (A' \bwt B')$. Cf. Appendix \cite{nrm} it follws that $L=\bigoplus_{dij \in \ca_{L}} B^{d}_{ij}\bwt ka^{d}_{ij}$. Then $f^{d}_{ij}\bwt a^{d}_{ij}\in A\bwt B$ implies that $f^{d}_{ij}\in A\cap A'$.\epf
In this way I can get Lagrangian subcategories $\cd(K,K)'=\cd(K,K)$.

\subsection{Other questions}

1. Normal coideal subalgebras of $D(A)$.

For $D(G)$ one does not have to ask for normality of the groups.

2. Quotients has $(N)$ would imply Hopf subalgebra have $(N)$ and viceversa.

3. Try to do as much as from the paper with $D(G)$.

4. Disccusions about orders of element and the decomposition in prime powers from \cite{DGNO}.

5. Need to think about a real application from dgno or weakly internal otherwise the paper is not very strong.

6. Definition of $ZLKer_{A}$ and its properties.

7. Is $$Z_D(M)\cong Z_A(M)\bowtie Z_A^*(M)?$$ Take $M$ to be simple first.

8. \bt\label{dualdefcl} 
Let $A$ be a modular Hopf algebra. Decompositions $\Rep(A)=\Dc \boxtimes \Dc'$ are in bijection with exact Hopf factorizations $A^*\cong (A//L)^*\otimes (A//C(L))^*$.
\et

9. In what conditions $D(N,\;M)$ is a normal coideal subalgebra. It is coideal but to be subalgebra requires few other things.
{\bf The proof is incomplete}
\bpf
Suppose as above that $\Dc=\Rep(A//L)$. Suppose that $\Dc'=\Rep(A//C(L))$. Then the multiplication map $(A//L)^*\otimes (A//C(L))^* \ra A^*$ is an isomorphism of coalgebras, therefore an exact factorization.

Conversely, suppose that $(A//L)^*\otimes (A//C(L))^* \ra A^*$ is an exact factorization. 
\epf

9. {\bf Another example}
$K=B(M,\;N,\;\lam)$ and $L=B(A,\;B,\;\mu)$. When these two commute?

10. $C_{\ch}\bwt C_d$ is closed under the adjoint action if and only if $\ch \bwt \mu$ is central.
\subsection{Aknowledgments}
This work was supported by a grant of the Romanian National Authority for Scientific Research, CNCS Ð UEFISCDI, project number PN-II-RU-TE-2012-3-0168.
\bibliographystyle{plain}
\bibliography{mec-rep-drd}
\ed

\bibliographystyle{plain}
\bibliography{mac-bob}
\end{document}